\documentclass[12pt]{article}

\usepackage{amsmath, amsthm, amssymb, amsfonts, enumerate}

 \usepackage{xcolor}

\def \R{I\!\!R}
\def \E{I\!\!E}

\def \N{I\!\!N}

\newtheorem{thm}{Theorem}[section]
\newtheorem{cor}[thm]{Corollary}
\newtheorem{lem}[thm]{Lemma}
\newtheorem{pro}[thm]{Proposition}
\newtheorem{defi}[thm]{Definition}
\newtheorem{rem}[thm]{Remark}

\newtheorem{rems}[thm]{Remarks}
\newtheorem{exa}[thm]{Example}

\def \conv{ {\rm conv} }

\def \conv{ {\rm conv} }

\def \Max{ {\rm Max} }

\def \lim{ {\rm lim} }

\def \argmax{ {\rm Argmax} }
\def \argmin{ {\rm Argmin} }
\def \Supp{ {\rm Supp} }

\def\abstract{\begin{center} \small\bf Abstract\end{center}\small}

 \title{Acyclic Gambling Games}
 
\author{Rida Laraki\thanks{LAMSADE, CNRS, [UMR 7243], Université Paris-Dauphine, PSL Research University, 75016 Paris, France . Part-time Professor at Ecole Polytechnique (Economics Department), Paris France. Honorary Visiting Professor at Liverpool University (Management School). Email:
rida.laraki@lamsade.dauphine.fr. Laraki's work was supported by (Idex [Grant Agreement No. ANR-11-IDEX-0003-02/Labex ECODEC No. ANR—11- LABEX-0047] and ANR-14-CE24-0007-01 “CoCoRICo-CoDec”).}  \; and J\'er\^ome Renault\thanks{Toulouse School of Economics, University of Toulouse Capitole, Toulouse, France. Email: jerome.renault@tse-fr.eu.}}

\date{\today} 

\begin{document}     

\maketitle
\begin{abstract}  

We consider   $2$-player zero-sum stochastic games where each player controls his own state variable living in a compact metric space. The terminology comes from gambling problems where the state of a player represents its wealth in a casino. Under natural assumptions (such as continuous running payoff and non expansive transitions), we consider for each discount factor the value $v_\lambda$ of the $\lambda$-discounted stochastic game and investigate its limit when $\lambda$ goes to 0. We show that under a strong acyclicity condition, the limit exists and is characterized as the unique solution of a system of functional equations: the limit is the unique continuous excessive and depressive function such that each player, if his opponent does not move, can reach the zone when the current payoff is at least as good than the limit value, without degrading the limit value. The approach generalizes and provides a new viewpoint on the Mertens-Zamir system coming from the study of zero-sum repeated games with lack of information on both sides. A counterexample shows that under a slightly weaker notion of  acyclicity, convergence of $(v_\lambda)$ may fail.

%\noindent {\it Key words.}  

\end{abstract}

\section{Introduction}

The model of zero-sum stochastic games was introduced by Shapley \cite{Sh} in 1953. A state variable $\omega\in \Omega$ follows a controlled Markov chain with transitions  $Q(\tilde \omega | i,j, \omega)$   controlled by the actions of two competing players ($i\in I$ for player 1 and $j\in J$ for player 2).  Shapley assumed the action and state spaces ($I$, $J$ and $\Omega$) to be finite and proved the existence of the value $v_{\lambda}$ of the $\lambda$-discounted game using a dynamic programming principle,  and characterized $v_\lambda$ as the unique fixed point of what has been called  {\it the Shapley operator} \cite{RS_2001, Sorin_2002}. 

Bewley and Kohlberg \cite{BK}, using algebraic tools, proved the existence of the asymptotic value $v=\lim_{\lambda \rightarrow 0} v_{\lambda}$. Actually, when action and state spaces are finite, the equations that define $v_{\lambda}$ may be described by finitely many polynomial inequalities, implying that $v_{\lambda}$ is semi-algebraic and so is converging. The extension of this result to infinite stochastic games is a central question in mathematical game theory   \cite{LS_2014, MertensSorinZamir_2015, NeymanSorin_2003, Sorin_2002}. 

Recently, several important  conjectures \cite{Mertens_1987, MertensSorinZamir_2015} were proved to be false. Vigeral \cite{Vigeral_2013} and Ziliotto \cite{Ziliotto_2016} provided   examples where the family $\{v_{\lambda}\}$ diverges as $\lambda$ approaches zero.  In Vigeral,  the state space $\Omega$ is finite  and the action sets $I$ and $J$ are  semi-algebraic. In Ziliotto, the set of actions is finite but the state space $\Omega$ is compact, and can be seen  as the space of common beliefs on a finite  state variable, controlled but not observed by the players. 

%The common feature between the examples is the oscillation of the transitions under optimal play \cite{SorinVigeral_2015}. 

On the other hand, there are many classes of stochastic games with general state space and action sets where $\{v_{\lambda}\}$ converges (for a recent survey, see \cite{LS_2014}). Many have in common some  \textit{irreversibility} in the  transitions. In recursive games \cite{Everett, SorinVigeral_2013} the current payoff is zero until the game is absorbed. In absorbing games \cite{Kohlberg_1974, MNR_2009, RS_2001} there is only one non-absorbing state. In repeated games with  incomplete information  \cite{AM_1995, MertensZamir_1971, RS_2001}, once a player reveals some information, he cannot withdraw it.  Similarly in splitting games \cite{Laraki_2001a, Laraki_2001b, Sorin_2002} the state follows a martingale which eventually converges. Interestingly, in all those classes of ``irreversible'' stochastic games, not only we have convergence but also an explicit characterization of the asymptotic value.  This leads many to anticipate that irreversibility has to do with convergence.

% (see for instance \cite{LS_2014}, page 80).

%This leads community to expect that convergence holds in general only under rreversibility. 

%may oscillate infinitely often between two states. More precisely, the examples may be reduced to a game where each player controls the exit time from a state in his disfavor. His optimal exit moment (an integer variable) depends on the discount factor (a continuous variable). This creates jumps and oscillations as the discount factor goes to zero.

Our paper provides a weak and a strong definition of irreversibility (we call acyclicity) and prove that they constitute the frontier between convergence and divergence of $\{v_{\lambda}\}$: strong acyclicity guarantee convergence while the closely related weak acyclicity do not. To do so, we restrict ourself to a new class which embeds any product stochastic game \cite{FSV_2008} and naturally extends gambling houses. 
 
 A classical gambling house problem \cite{DS_1965, MS_1996} has three
ingredients :  a metric \textit{state space} $S$, a  Borel-measurable \textit{utility function} $u:S \rightarrow \R $, and a 
\textit{gambling house} $\Phi $, where $\Phi $ is a set value function that assigns to each $s\in S$ a set $\Phi (s)$ of $\Delta(S)$ (set of Borel probability distributions over $S$). At each stage $t$, given the state $s_t$, the decision maker gets the reward $u(s_t)$, chooses $p_t \in \Phi(s_t)$, and the state moves to $s_{t+1}$ according to the probability distribution $p_t$. 

The gambling house is \textit{leavable} if for every $s\in S$, $\delta_s$ (the dirac mass at $s$) belongs to $\Phi (s)$. It is well known that any MDP can be mapped to a gambling house and any positive MDP to a leavable gambling house \cite{DMS, MS_1996, Schal}.   In this paper, we consider only leavable gambling houses.

In a gambling {\it  game}, each player controls his gambling house: $\Gamma : X \rightarrow \Delta(X)$ for Player 1 and $\Lambda : Y \rightarrow \Delta(Y)$ for Player 2, and the utility function is now $u:X \times Y \rightarrow \R$ with the convention that player 1 wants to maximize $u$ whereas player 2 wants to minimize $u$. At each stage $t$, both players knowing the state $\omega_t=(x_t,y_t)$,  simultaneously Player 1 chooses $p_t$ in $\Gamma(x_t)$ and  Player 2 chooses $q_t \in \Lambda (y_t)$, the stage payoff is $u(x_t,y_t)$ and a new state $(x_{t+1},y_{t+1})$ is selected according to the probability distribution $p_t \otimes q_t$. As any MDP can be mapped into a gambling house, a product stochastic game \cite{FSV_2008} can be  mapped  into a gambling game. 

For each $\lambda \in (0,1]$,  one can define the $\lambda$-discounted game where the stream of payoffs is evaluated according to $\sum_{t=1} \lambda(1- \lambda)^{t-1} u(x_t,y_t)$. Usually $\lambda$ is called the discount  rate, $1-\lambda=\frac{1}{1+r}$ is called the discount factor  and $r$ is the interest rate. Under classical regularity assumptions, the $\lambda$-discounted game has a value $v_{\lambda}$  and the family $\{ v_{\lambda}\}$ is equi-continuous.\\

%\footnote{A recent Tauberian theorem of Ziliotto \cite{Ziliotto_2015} shows that, under equi-continuity, convergence/divergence of $v_{\lambda}$ extends to $v_n$ where $v_n$ is the value of $n$-stage gambling game where the stream of payoffs is evaluated as $\frac{1}{n} \sum_{t=1}^n u(x_t,y_t)$. More precisely, $\lim_{n\rightarrow \infty}v_n$  exists if and only if  $\lim_{\lambda \rightarrow 0}v_{\lambda}$ exists and whenever they exist, they coincide.}. 

Our first main result shows that if at least one of the gambling houses   $\Gamma$ or $\Lambda$ is strongly acyclic, $\{v_{\lambda}\}$ uniformly converges to a function $v$ as $\lambda$ goes to $0$. Moreover, we provide several characterizations of the asymptotic value $v$ that extend the well known Mertens-Zamir system of functional equations \cite{MertensZamir_1971}.  Our second result proves that under a slightly weaker notion of acyclicity, $\{v_{\lambda}\}$ may diverge (even if both houses $\Gamma$ and $\Lambda$ are weakly acyclic and both state spaces $X$ and $Y$ are finite). Our example is inspired by an example in Ziliotto \cite{Ziliotto_2016} for stochastic games where both players control the same state variable. It is the first in the class of product stochastic games, and appears  simpler than the existing counterexamples of divergence with finite state space.  Finally, under an idempotent assumption combined with a bounded variation hypothesis on the transitions, we prove existence of the uniform value, extending a recent result by Oliu-Barton \cite{Miquel} on splitting games.

\section{Gambling Games}

\subsection{Notations}

Given a compact metric space $S$, we denote by $\mathcal{B}(S)$, resp. by $\mathcal{C}(S)$,   the  set of bounded measurable, resp. continuous, functions from $S$ to the reals, and by $\Delta(S)$     the set of Borel probabilities over $S$.  For $s$ in $S$,  we denote by  $\delta_s\in \Delta(S)$ the Dirac measure on $s$, and whenever  possible we  assimilate $s$ and $\delta_s$. For $v$ in $B(S)$, we denote by $\tilde{v}$ its affine extension to  $\Delta(S)$: $\tilde{v}(p)=\E_p(v)$ for all $p$ in $\Delta(S)$, where $\E_p(v)$ is the expectation of $v$ with respect to $p$.  $\Delta(S)$ is endowed with the weak-* topology, a compatible distance  being  the Kantorovich-Rubinstein (or Wasserstein of order 1) metric: $d_{KR}(p,p')= \sup_{v\in E_1} |\tilde{v}(p)-\tilde{v}(p')|$, where $E_1$ is the set of 1-Lipschitz functions on $S$.  When there is no confusion, $\tilde{v}(p)$ will also be denoted by $v(p)$.

\subsection{Model}

A gambling game is a zero-sum stochastic game where each player controls  his own state variable. We will always assume in this paper that the state spaces are non empty metric and compact, and   denote by  $X$ and  $Y$ the respective set of states\footnote{Both metrics are denoted by $d$, and we will use the metric $d((x,y),(x',y'))=d(x,x')+d(y,y')$ on $X\times Y$.}  controlled by Player 1 and by Player 2. The  transitions of Player 1 are given by a continuous\footnote{i.e. $\forall \varepsilon>0, \exists \alpha>0, \forall x, x' \in X \; {\rm with}\;  d(x,x')\leq \alpha, \forall p \in \Gamma(x),  \exists p'\in \Gamma(x') \; {\rm s.t.}\;  d_{KR}(p,p')\leq \varepsilon$. Since $\Gamma$ has compact values, it implies that the graph of $\Gamma$ is compact.} multifunction $\Gamma: X \rightrightarrows \Delta(X)$ with non empty convex\footnote{If originally the gambling game has non  convex values, then allowing as usual  players to randomize, would lead to transitions with convex values.} 
  compact values: if the state of Player 1 is at $x$,   he can select his new  state according to any probability in $\Gamma(x)$. Similarly, a continuous multifunction $\Lambda:  Y \rightrightarrows \Delta(Y)$ with non empty convex compact values, gives the transitions of Player 2. The players independently control their own state, and  only interact through payoffs: the running payoff of Player 1 is given  by a continuous mapping $u:X\times Y \longrightarrow \R$, and the payoff to Player 2 is given by $-u$.

Gambling games extend the model of gambling houses \cite{DS_1965}, which correspond to the single player case  when $Y$ is a singleton and Player 2 plays no role. 
One can show, by an adequate increase of the state space in order to   encompass  actions,   that any MDP can be mapped into   a gambling house \cite{DMS, MS_1996, Schal}. 

A standard gambling house is the red-and-black casino where $X=[0,1]$ is a  fortune space. Suppose that at each fortune $x\geq 0$, the gambler
can stake any amount $s$ in her possession. The gambler loses the stake with
probability $1-w$ where $w\in (0,1)$ is fixed and given,  and wins back the stake and an additional equal amount   with probability $w$. The corresponding  transition multiifunction reads: 
\begin{equation*}
\Phi_{w}(x)=\{w\delta_{\min\{x+s,1\}}+(1-w)\delta_{x-s}\,:\,0\leq s\leq x\}.
\end{equation*}
Another class  of gambling house are splitting problems where $X=\Delta(K)$ is a simplex ($K$ is a finite set) and $\Gamma(x)$ is the set of Borel probabilities on $X$ centered at $x$. The idea of splitting was introduced by Aumann and Maschler \cite{AM_1995} in the context of repeated games with incomplete information on one side. This is now very popular in economics (persuasion and information design literature \cite{KG}).

The above examples of gambling houses naturally extend to gambling games. 
One can consider a casino game where each player $i$ controls a red-and-black house with parameter $w_i$, and the running payoff depends on the current pair of fortunes. Another example is a splitting game  \cite{Laraki_2001a, Laraki_2001b, Sorin_2002} where $X=\Delta(K)$ and $Y=\Delta(L)$ are simplexes, 
 $\Gamma(x)$ is the set of Borel probability measures on $X$ that are centered at $x$ and  $\Lambda(y)$ is the set of probability measures of Borel probability measures on $Y$ that are centered at $y$.  
 
 %These splitting games comes from repeated games with lack of information on both sides (Aumann Maschler, Mertens Zamir) :  here $x\in \Delta(K)$ is controlled by Player 1 and  represents the uncertainty of Player 2 about some fixed state of nature $k\in K$ known by Player 1, and $\Gamma(x)$ represents the set of beliefs that Player 1 can induce by sending  messages to Player 2. Similarly $y\in \Delta(L)$   represents the uncertainty of Player 1 about some state of nature $l$ known by Player 2, and $\Lambda(y)$ represents the set of beliefs that Player 2 can induce by sending  messages to Player 1. 

\subsection{Discounted Evaluations}

Given a discount factor $\lambda\in (0,1]$ and an initial state  $(x_1,y_1)$ in $X\times Y$, the game $G_\lambda(x_1,y_1)$ is played as follows: at any stage $t\geq 1$, the payoff to Player 1   is $u(x_t,y_t)$ and both players knowing $(x_t,y_t)$, simultaneously Player 1  chooses $p_{t+1}$ in $\Gamma(x_{t})$ and Player 2 chooses $q_{t+1}$ in $\Lambda(y_{t})$. Then,   $x_{t+1}$ and $y_{t+1}$ are independently selected according to $p_t$ and $q_t$,  the new states $x_{t+1}$ and $y_{t+1}$ are publicly announced, and the play goes to stage $t+1$. Under our assumptions of compact state  spaces, continuous transitions with convex compact values and continuous running payoff, it  is easy to describe  the value of such dynamic game.

\begin{defi}  \label{def0}  $v_\lambda$ is  the unique   element of  $\mathcal{C}(X\times Y)$ s.t.  $\forall (x,y)\in X\times Y,$
\begin{eqnarray*}
v_\lambda(x,y)&=&\max_{p \in \Gamma(x)}\min_{q \in \Lambda(y)} \left(\lambda u(x,y)+ (1-\lambda) \tilde{v}_\lambda(p,q)\right),\\
&=&\min_{q \in \Lambda(y)}\max_{p \in \Gamma(x)} \left(\lambda u(x,y)+ (1-\lambda) \tilde{v}_\lambda(p,q)\right).
\end{eqnarray*}

 \end{defi}
\noindent This is the standard characterization of the value of a discounted game by means of the Shapley operator. Existence and uniqueness of $v_{\lambda}$ follow from standard fixed-point arguments (see for instance \cite{MertensSorinZamir_2015, RS_2001}). We refer to $v_\lambda(x,y)$ as  the value of the game  $G_\lambda(x,y)$.\\

The goal of the paper is to study the convergence of $(v_\lambda)_\lambda$ when $\lambda$ goes to 0, i.e. when players become more and more patient. 

\begin{rem}\rm Cesaro Evaluations. It is also standard to define the value of the $n$-stage games by: $v_1=u$, and for $n\geq 1$ and  $(x,y)\in X\times Y$:\begin{eqnarray*}
v_{n+1}(x,y)&=&\frac{1}{n+1} \max_{p \in \Gamma(x)}\min_{q \in \Lambda(y)} \left(  u(x,y)+ n  \tilde{v}_n(p,q)\right),\\
&=&\frac{1}{n+1} \min_{q \in \Lambda(y)}\max_{p \in \Gamma(x)} \left( u(x,y)+n \tilde{v}_n(p,q)\right).
\end{eqnarray*}
It is known that the uniform convergence of $(v_n)_n$  when $n$ goes to infinity, is equivalent  to the uniform convergence of $(v_\lambda)_\lambda$ when $\lambda$ goes to 0, and in case of convergence both limits are the same (Theorem 2.2 in \cite{Ziliotto_2015} applies here)\end{rem}

\subsection{Non expansive transitions}
 Without further assumptions, convergence of $(v_\lambda)$ may fail even in the simple case where $\Gamma$ and $\Lambda$ are single-valued (``0 player case", players have no choice), so we will assume throughout the paper that the gambling game is non expansive, i.e. have non expansive transitions: 

 \begin{defi} \label{def0.5} The game has non expansive transitions if: 
$$\forall x\in X, \forall x'\in X, \forall p \in \Gamma(x), \exists p'\in \Gamma(x'),\; s.t. \;  d_{KR}(p,p')\leq d(x,x'),$$
   $${\it and\; similarly:}\;  \forall y\in Y, \forall y'\in Y, \forall q \in \Lambda(y), \exists q'\in \Lambda(y'),\; s.t. \;  d_{KR}(q,q')\leq d(y,y').$$
   \end{defi}
   
The gambling game has non expansive transitions if the transitions, viewed as mappings from $X$ to $2^{\Delta(X)}$, and from  $Y$ to $2^{\Delta(Y)}$,  are 1-Lipschitz for the Hausdorff distance on compact subsets of ${\Delta(X)}$ and ${\Delta(Y)}$.  Note that the transitions are always non expansive when $X$ and $Y$ are finite\footnote{If  $X$ is finite and $d(x,x')=2$ for $x\neq x'$, then $d_{KR}(p,p')=\|p-p'\|_1$ for $p$, $p'$ in $\Delta(X)$.}. Moreover splitting games are non expansive \cite{Laraki_2001a},  and red-and-black casino houses with parameter $w$ are non-expansive if and only if $w\leq \frac{1}{2}$ \cite{LS_2004}.  Let us mention also  Markov chain repeated games  with incomplete information \cite{GR_2015}, where each player observes a private and  exogenous   Markov chain. These repeated games  lead to gambling houses with   transitions of the form: $X$ is a simplex  $\Delta(K)$, and $\Gamma(x)=\{pM, p \in \Delta(X) \; {\rm centered \; at \;} x\}$ with $M$ a fixed stochastic matrix. Here again, transitions are non expansive.\\ 

Let us mention immediately an important consequence of the non expansive assumption. The proof is in the Appendix. 

 \begin{pro} \label{pro1} Assume the gambling game is non-expansive. Then the family $\{v_\lambda\}_{\lambda\in (0,1]}$ is equicontinuous. \end{pro}
 
 This proposition extends to two players a similar result in \cite{LS_2004} on gambling houses where it is proved that non-expansivity is necessary and sufficient to guarantee equi-continuity of the values. As a consequence, pointwise and uniform convergence of $\{v_\lambda\}$ are equivalent, and since $X\times Y$ is compact, to prove this convergence it is enough to prove uniqueness of a limit point\footnote{By convergence or limit point of  $\{v_\lambda\}$, we always mean when $\lambda$ approaches 0.}.
 
\begin{rem}\rm It is not difficult to see that without non-expansivity, $\{v_\lambda\}$ may not be equicontinuous and the convergence   may not be uniform.   For instance  in red-and-black casino with a single player, if the parameter $w>\frac{1}{2}$   and $u(x)=x$, $v_{\lambda}$ is continuous for every $\lambda$ but $v=\lim_{_\lambda \rightarrow 0} v_{\lambda}$ is not  : $v(x)=0$ for $x=0$ and $v(x)=1$ for $x>0$.\end{rem}

 \subsection{Excessive, depressive and balanced functions}
 
 \begin{defi} \label{def2}$\;$ Let $v$ be in ${\cal B} (X\times Y)$.
 
 1) $v$ is balanced if $\forall (x,y)\in X \times Y$, 

\centerline{$v(x,y)=\max_{p\in \Gamma(x)} \min_{q \in \Gamma(y)} \tilde{v}(p,q)=  \min_{q \in \Gamma(y)}\max_{p\in \Gamma(x)} \tilde{v}(p,q).$}

2) $v$ is excessive (with respect to $X$) if: $\forall (x,y)\in X \times Y$, $v(x,y)=\max_{p \in \Gamma(x)} \tilde{v}(p, y).$

3) $v$ is depressive (with respect to $Y$)  if: $\forall (x,y)\in X \times Y$, $v(x,y)=\min_{q \in \Gamma(y)} \tilde{v}( x,q).$

\end{defi}

Observe that any uniform limit $v$ of  $(v_\lambda)_{\lambda\in (0,1]}$ is necessarily continuous and balanced (by passing to the limit in definition \ref{def0}).  
In a splitting game, excessive means concave with respect to the first variable, and depressive means convex with respect to the second variable.   

\begin{exa} \label{exa3} \rm Consider a gambling game where players 1 and 2 move  on the same finite grid of a circle containing $6$ nodes in equidistant positions. Any player can move one step to the left, one step to the right, or not move (and choose randomly between these 3 options, so that transitions have convex values).  It is here possible for a player to go from any state to any other state in at most $3$ stages (the game may be called \textit{cyclic}),  so any excessive and depressive function is necessarily constant. Suppose that Player 1's payoff is 1 if he is at  most one step away from Player 2, and Player 1's payoff is   0 otherwise. If the players start a distance at most 1, Player 1 can guarantee this property will hold forever   by not moving or moving one step to the direction of Player 2 and so, in this case we have $v_{\lambda}=1$ for every $\lambda$. On the other hand, if the players start at a distance at least  $2$, Player 2 can insure that this property will hold forever, by not moving or moving one step  in the opposite direction of Player 1. For these initial states, $v_{\lambda}=0$  for every $\lambda$. Here,  $v=\lim \, v_{\lambda}$ is continuous and balanced, but not  excessive nor depressive.  
 \end{exa}

\subsection{Acyclicity}

We now come to the main definitions of the paper.

%Recall that in splitting games (Sorin 2002, Laraki 2001 and 2002), $\{v_\lambda\}$ uniformly converges to the unique continuous concave-convex function $v:X\times Y\longrightarrow \R$ satisfying the Mertens-Zamir (1971, 1976) system of functional equations:
% $$MZ1: v= Cav_{X} \min (u,v), $$
%$$MZ2: v= Vex_{Y} \max (u,v),$$
%where $Cav_{X}$ is the concavification operator on $X$ which associates to every function $f \in \mathcal{B}(X)$ the smallest concave function larger than $f$, and similarly, $Vex_{Y}$ is the convexification operator on $Y$ which associates to every function $g \in \mathcal{B}(Y)$ the highest convex function smaller than $g$. 
%
%Under acyclicity and non-expansiveness assumptions on $\Gamma$ and $\Lambda$, we will prove uniform convergence of $v_\lambda$ (when $\lambda$ goes to 0) to a limit value function $v$ (called the asymptotic value), and  give a new variational characterization for $v$, that extends $MZ$. Moreover, we will show that our result is tight in the sense that a slight weakening of acyclicity leads to divergence of $v_{\lambda}$. 

\begin{defi} \label{def1}$\;$

1) The gambling game is leavable if:  $\forall x\in X,\delta_x \in \Gamma(x)\; {\it  and} \;   \forall y\in Y, \delta_y \in \Lambda(y).$

2) The gambling house $\Gamma$ of player 1 is weakly acyclic if there exists $\varphi$ in ${\cal B}(X)$  lower semi-continuous  such that:
$$\forall x\in X, \argmax_{p \in \Gamma(x)} \tilde{\varphi}(p)=\{\delta_x\}.$$
Similarly, the gambling house $\Lambda$ of player 2 is weakly acyclic  if there exists  $\psi$ in ${\cal B}(Y)$  upper semi-continuous   such that:
$$\forall y\in Y, \argmin_{q \in \Lambda(y)} \tilde{\psi}(q)=\{\delta_y\}.$$
The gambling game is weakly acyclic if both gambling houses are weakly acyclic. \end{defi}

The gambling game is leavable if each player can remain in any given state. This is a standard assumption \cite{DS_1965}. This is the case in red-and-black casinos and splitting games. In persuasion games and in repeated games with incomplete information, not moving means  revealing no information. If the game is leavable, any excessive and depressive function is necessarily balanced. The converse is not necessarily true as example \ref{exa3} shows.

Weak acyclicity is, to our knowledge, a new condition in the gambling house literature. If  the house $\Gamma$ is weakly acyclic,  the ``potential" $\varphi$ decreases  in expectation along  non stationary trajectories, hence the irreversibility of the process (in the space of probabilities over $X$). 

%Thus , it acts a Lyapunov function, implying convergence of every trajectory.  

Observe that any weakly acyclic gambling game is necessarily leavable. When $w\leq \frac{1}{2}$, a red-and-black casino is weakly acyclic (take $\varphi$ to be   strictly increasing and strictly concave). Also, a splitting game is weakly acyclic (take $\varphi$ to be any strictly concave function on $X$).\\

We will now define strong acyclicity, our central condition. For this, we need to consider transitions for several stages. We first extend linearly the transitions to $\Delta(X)$ and $\Delta(Y)$ by defining $\tilde{\Gamma}: \Delta(X) \rightrightarrows \Delta(X)$ and $\tilde{\Lambda}: \Delta(Y) \rightrightarrows \Delta(Y)$. 
More precisely, the graph of  $\tilde{\Gamma}$ is defined as the closure of the  convex hull of the graph of $\Gamma$ (viewed as the subset $\{(\delta_x,p), x\in X, p \in \Gamma(x)\}$ of $\Delta(X)\times \Delta(X)$),  and  similarly the graph of  $\tilde{\Lambda}$ is defined as the closed  convex hull of the graph of $\Lambda$. Because Dirac measures are extreme points of $\Delta(X)$ and $\Delta(Y)$, we have $\tilde{\Gamma}(\delta_x)=\Gamma(x)$ and $\tilde{\Lambda}(\delta_y)=\Lambda(y)$ for each $x$ in $X$  and $y$ in $Y$.  
Be careful that in general, for $p$ in $\Delta(X)$ and $q$ in $\Delta(Y)$: $\tilde{v}_\lambda(p,q) \neq  \max_{p' \in \tilde{\Gamma}(p)}\min_{q' \in \tilde\Lambda(q)} \left(\lambda u(p,q)+ (1-\lambda)  {\tilde {v}}_\lambda(p',q')\right).$\\

We  now define inductively a sequence of transitions $(\tilde{\Gamma}^n)_n$ from $\Delta(X)$ to $\Delta(X)$, by   $\tilde{\Gamma}^0(p)=\{p\}$ for every state $p$ in $\Delta(X)$, and\footnote{The composition being  defined by $G\circ H(p)=\{p"  \in\Delta(X), \exists p'\in H(p)\; {\rm s.t.}\;  p"\in G(p')\}$.}  for each $n\geq 0$, $\tilde{\Gamma}^{n+1}=\tilde{\Gamma}^n\circ \tilde{\Gamma}$ . $\tilde{\Gamma}^n(\delta_x)$  represents the set of probabilities over states that Player 1 can reach in $n$ stages from the initial state  $x$ in $X$. Similarly we define $\tilde{\Lambda}^n$ for each $n$.

 \begin{defi}  \label{def3}  $\;$
 
 1) The reachable set of Player 1 from state $x$ in $X$  is  the closure of  $\bigcup_{n\geq 0} \tilde{\Gamma}^n(\delta_x)$ in $\Delta(X)$, and denoted $\Gamma^\infty(x)$.
  Similarly, the reachable set of Player 2 from state $y$ in $Y$ is the subset $\Lambda^\infty(y)$ of $\Delta(Y)$ defined as the closure of  $\bigcup_{n\geq 0} \tilde{\Lambda}^n(\delta_y)$.\\

2)  The gambling house $\Gamma$ of player 1 is strongly acyclic (or simply, acyclic) if there exists $\varphi$ in ${\cal B}(X)$ lower semi-continuous  such that:
$$\forall x\in X, \argmax_{p \in \Gamma^\infty(x)} \tilde{\varphi}(p)=\{\delta_x\}.$$
Similarly, the gambling house $\Lambda$ of player 2 is strongly acyclic (or simply, acyclic)   if there exists  $\psi$ in in ${\cal B}(Y)$ upper semi-continuous   such that:
$$\forall y\in Y, \argmin_{q \in \Lambda^\infty(y)} \tilde{\psi}(q)=\{\delta_y\}.$$
 The  gambling game is   strongly acyclic (or simply, acyclic) if both gambling houses are strongly acyclic.

 \end{defi}

Clearly, strong  acyclicity implies weak  acyclicity. Note that in splitting games, we have  $\Gamma \circ \Gamma = \Gamma$,  so $\Gamma^\infty=\Gamma$  and weak and strong acyclicity coincide.

 \section{Main Results}

In all the paper (and in particular the results  below), we consider {\it standard}  gambling games.  

\begin{defi}
A gambling game is standard if both state spaces $X$ and $Y$ are compact metric, the running payoff  $u$ is continuous, and the transitions $\Gamma$ and $\Lambda$  have non empty convex compact values and are leavable and non expansive.
\end{defi}

We will also   use the following properties.  

\begin{defi} Given $v$ in $\mathcal{B}(X\times Y)$, we say that :

$v$ satisfies  $P1$ if:  $\forall (x,y)\in X\times   Y, \exists p \in \Gamma^\infty(x), v(x,y) = \tilde{v}(p,y)\leq u(p,y),$

$v$ satisfies $P2$ if : $\forall (x,y)\in X\times   Y, \exists q \in \Lambda^\infty(y), v(x,y)= \tilde{v}(x,q)\geq u(x,q).$

\end{defi}

% \begin{defi}  \label{defx} \ \\
% The game is one-side (resp. both-sides) strongly acyclic if only one of the houses is (resp. both houses are) strongly acyclic.\ \\
% The game is one-side (resp. both-sides) weakly acyclic if only one of the houses is (resp. both houses are) weakly acyclic.
% \end{defi}

 Our main result is the following.

\begin{thm} \label{thm1} Consider a standard gambling game.

\noindent 1. If at least one of the players has a  strongly acyclic gambling house,  $(v_\lambda)$ uniformly converges to the unique function $v$ in $\mathcal{C}(X\times Y)$ satisfying: \\

a) $v$ is excessive , i.e. $\forall (x,y)\in X \times Y$, $v(x,y)=\max_{p \in \Gamma(x)} \tilde{v}(p, y),$\\

  b) $v$ is depressive, i.e.  $\forall (x,y)\in X \times Y$, $v(x,y)=\min_{q \in \Gamma(y)} \tilde{v}( x,q),$\\

 c) $v$ satisfies $P1$, i.e.  $\forall (x,y)\in X\times   Y, \exists p \in \Gamma^\infty(x), v(x,y) = \tilde{v}(p,y)\leq u(p,y),$\\
 
d) $v$ satisfy $P2$, i.e.  $\forall (x,y)\in X\times   Y, \exists q \in \Lambda^\infty(y), v(x,y)= \tilde{v}(x,q)\geq u(x,q).$\\

\noindent Moreover: $v$ is the largest excessive-depressive continuous function satisfying $P1$, and is the smallest excessive-depressive continuous function satisfying $P2$.\\
 
\noindent 2.  If both  gambling houses are weakly acyclic, convergence of $(v_\lambda)$ may fail.

\end{thm}

\vspace{1cm}

The   conditions of the positive result 1) may be interpreted as follows:

\begin{itemize}
\item a) and b) : {\it It  is always safe not to move.} For each player, not moving ensures not to degrade the limit value. 

\item c) and d) :  {\it Each player can reach, if his opponent does not move,   the zone when the  current payoff is at least as good  than the limit value, without degrading the limit value.}

\end{itemize}

\noindent These interpretations will lead later  to the  construction of simple uniformly optimal strategies in  some gambling games, see section \ref{subuniform}. \\

The  positive result of theorem \ref{thm1}   relies on  the following   three propositions. Recall that thanks to proposition \ref{pro1},  to get convergence of the values it is enough to show uniqueness of a limit point of $(v_\lambda)_\lambda$.

 \begin{pro}  \label{pro2} Assume one of the player has a weakly acyclic gambling house. If $v$ in $\mathcal{C}(X\times Y)$ is balanced, then $v$ is excessive and depressive.  \end{pro}

 \begin{pro} \label{pro2.5} Let $v$ be a limit point of $(v_\lambda)$ for the uniform convergence. Then $v$ is balanced, and satisfies $P1$ and $P2$. \end{pro}
 
% We have:   
% $$v \; {\it is} \;  {\it balanced},$$
% $$P1: \forall (x,y)\in X\times   Y, \exists p \in \Gamma^\infty(x), v(x,y)= v(p,y)\leq u(p,y),$$
%$$P2: \forall (x,y)\in X\times   Y, \exists q \in \Lambda^\infty(y), v(x,y)= v(x,q)\geq u(x,q).$$ \end{pro}

\begin{pro} \label{pro3} Assume  one of the player has a strongly  acyclic gambling house. Then,  any balanced continuous function satisfying  $P1$ is smaller that any balanced continuous function satisfying  $P2$. Consequently, there is at most one balanced continuous function satisfying  $P1$ and $P2$. 

%$$P1: \forall (x,y)\in X\times   Y, \exists p \in \Gamma^\infty(x), v(x,y) \leq v(p,y)\leq u(p,y),$$
%$$P2: \forall (x,y)\in X\times   Y, \exists q \in \Lambda^\infty(y), v(x,y) \geq v(x,q)\geq u(x,q).$$
\end{pro}

\vspace{0.5cm}

If none of the player has a strongly  acyclic gambling house, there may be infinitely many balanced continuous functions satisfying P1 and P2. This will be the case in our counter-example of section  \ref{contrex}, where both gambling houses are weakly acyclic.

\section{Examples}

\subsection{A simple acyclic gambling house}\label{exa1}

Let us first  illustrate our characterization on a simple example.  Consider the following Markov decision process with 3 states: $X=\{a,b,c\}$ from \cite{Sorin_2002}.   States   $b$ and $c$ are absorbing with respective payoffs 1 and 0. Start at $a$, choose $\alpha\in I=[0,1/2]$, and move to $b$ with proba $\alpha$ and to $c$ with proba $\alpha^2$. 
 
 \vspace{0.5cm}
 
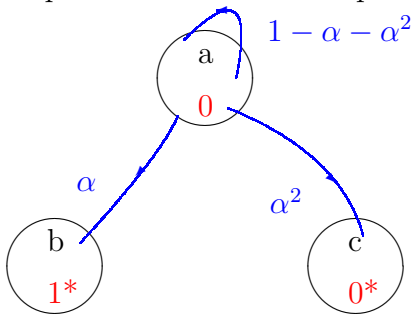
 \begin{center}
\setlength{\unitlength}{1mm}
\begin{picture}(60,30)

 \put(30,30){\circle{12}}
  \put(29,32){a}
   \textcolor{red} {    \put(29,25){0}}
    
     \put(10,5){\circle{12}}
  \put(9,7){b}
  \textcolor{red} {     \put(9,0){1*}}
    
         \put(50,5){\circle{12}}
  \put(49,7){c}
   \textcolor{red} {    \put(49,0){0*}}

    \textcolor{blue} {      { \qbezier(25,25)(23,20)(12,8)}}
 \textcolor{blue} {   { \put(19,18){\vector(-1,-1){1}}}}
 
     \textcolor{blue} {      { \qbezier(30,30)(33,45)(23,35)}}
 \textcolor{blue} {   { \put(27,39){\vector(-1,0){1}}}}
 
 \textcolor{blue} {     \put(6,15){$\alpha$}}
  \textcolor{blue} {     \put(30,35){$1-\alpha- \alpha^2$}}
   
    \textcolor{blue} {     \put(29,12){$\alpha^2$}}
 
 \textcolor{blue} {     { \qbezier(22,26)(37,20)(40,9)}}
 \textcolor{blue} {   { \put(34,17){\vector(1,-1){1}}}}

\end{picture}
    
      \end{center}

\noindent Here formally $Y$ is a singleton (there is only one player, so we can omit the variable $y$), $\Gamma(b)=\{\delta_b\}$, $\Gamma(c)=\{\delta_c\}$, and $\Gamma(a)=\conv\{(1-\alpha-\alpha^2) \delta_a+ \alpha \delta_b +\alpha^2 \delta_c, \alpha \in I\}.$ The payoffs are $u(a)=u(c)=0$, $u(b)=1$. 

$\Gamma$ has compact convex values, the transitions are 1-Lipschitz, and the game is leavable. The gambling game is strongly acyclic: just   consider $\varphi$ such that $\varphi(a)=1$, and  $\varphi(b)=\varphi(c)=0$.

Player 1 can   go from state $a$ to state $b$ in infinitely  many stages with arbitrarily high probability, by repeating a choice of $\alpha>0$ small (so that $\alpha^2$ is much smaller than $\alpha$), and the limit value $v$ clearly satisfies:
$$v(a)=v(b)=1, v(c)=0.$$

This is the unique function $w:X\rightarrow \R$ satisfying the  conditions a), b), c), d) of   Theorem \ref{thm1}:   $P1$ and $P2$ implies  $u\leq w\leq 1$, and because  $b$ and $c$ are absorbing states, $w(b)=1$ and $w(c)=0$.  Finally,  $w$ excessive gives $w(a)=1$.  Notice that $\delta_b \in \Gamma_\infty(a)$ but for each $n$, $\delta_b \notin \tilde{\Gamma}_n(a)$.

\subsection{A weakly acyclic gambling house} \label{exa2}

Let us modify the   gambling house of the previous section \ref{exa1}. We still have a unique player and 
a state space  $X=\{a,b,c\}$. The only difference is that     state $b$ is no longer absorbing : in state b the player  also has to choose some $\alpha \in I=[0,1/2]$, and  then moves to $a$ with probability  $\alpha$, to $c$ with probability  $\alpha^2$  and remains in $b$ with probability  $1-\alpha-\alpha^2$.

\vspace{1cm}

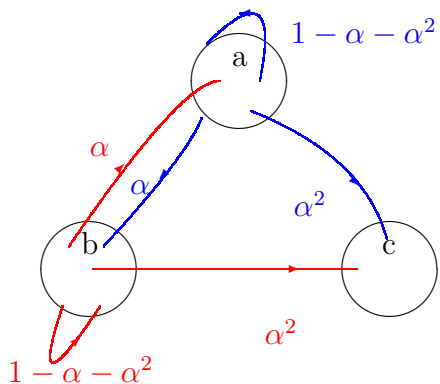
   \begin{center}
\setlength{\unitlength}{1mm}
\begin{picture}(60,30)

 \put(30,30){\circle{12}}
  \put(29,32){a}

     \put(10,5){\circle{12}}
  \put(9,7){b}

         \put(50,5){\circle{12}}
  \put(49,7){c}

    \textcolor{blue} {      { \qbezier(25,25)(23,20)(12,8)}}
 \textcolor{blue} {   { \put(19,18){\vector(-1,-1){1}}}}
 
     \textcolor{blue} {      { \qbezier(30,30)(33,45)(23,35)}}
 \textcolor{blue} {   { \put(27,39){\vector(-1,0){1}}}}
 
 \textcolor{blue} {     \put(10,15){$\alpha$}}
  \textcolor{blue} {     \put(30,35){$1-\alpha- \alpha^2$}}
   
    \textcolor{blue} {     \put(29,12){$\alpha^2$}}
 
 \textcolor{blue} {     { \qbezier(22,26)(37,20)(40,9)}}
 \textcolor{blue} {   { \put(34,17){\vector(1,-1){1}}}}

     \textcolor{red} {      { \qbezier(-5,8)(10,30)(15,30)}}
 \textcolor{red} {   { \put(0,18){\vector(1,1){1}}}}
  \textcolor{red} {     \put(-5,20){$\alpha$}}

     \textcolor{red} {      { \qbezier(-10,0)(-15,-15)(-5,0)}}
 \textcolor{red} {   { \put(-10,-5){\vector(1,1){1}}}}

  \textcolor{red} {     \put(-20,-10){$1-\alpha- \alpha^2$}}

 \textcolor{red} {     { \qbezier(-10,5)(15,5)(25,5)}}
 \textcolor{red} {   { \put(15,5){\vector(1,0){1}}}}
     \textcolor{red} {     \put(10,-5){$\alpha^2$}}
\end{picture}
    
      \end{center}

\vspace{1cm}

  States $a$ and $b$ are now symmetric. This gambling house  is weakly acyclic, with $\varphi(a)= \varphi(b)=1, \varphi(c)=0$, but  it is not strongly acyclic  since    $a\in \Gamma^\infty(b)$ and $b\in \Gamma^\infty(a)$. We will later use this gambling house to construct our counter-example of theorem \ref{thm1}, 2).

\subsection{An example with countable state spaces}

We present here a (strongly) acyclic gambling game with countable state spaces, and  illustrate\footnote{We have first studied and understood this  example before proving  proposition \ref{pro2}. Example \ref{exa3} shows that weakening the assumption of weakly acyclic gambling game to  leavable gambling game is not possible in this proposition.}  the proof of proposition \ref{pro2}, that under weak acyclicity  any continuous balanced function is also excessive and depressive. Consider the state space:
$$X= \{1-\frac{1}{n}, n \in \N^*\}\cup\{1\}=\{x_1,..., x_n,...., x_{\infty}\},$$
\noindent where $x_n=1-\frac{1}{n}$ if $n$ is finite, and $x_{\infty}=1$. We use $d(x,x')=|x-x'|$, so that $X$ is countable and compact. The transition is given by:
 $$\Gamma(x_n)=\{\alpha \delta_{x_n}+(1-\alpha)  \delta_{x_{n+1}}, \alpha \in [0,1]\}, \; {\rm and}\;  \Gamma(x_\infty)=\{\delta_{x_\infty}\}.$$
 The intuition is clear: Player 1 can stay at his location, or move 1 to the right. The gambling house $(Y,\Lambda)$ of 
  Player 2 is a copy of the gambling house  of Player 1. Transitions are non expansive (since $|\frac{1}{n+1}-\frac{1}{n'+1}|\leq |\frac{1}{n}-\frac{1}{n'}|)$, and the game is strongly acyclic. The payoff $u$ is any continuous function $X\times Y\longrightarrow \R$, so that  theorem \ref{thm1} applies.\\

Consider  $v:X\times Y\longrightarrow \R$, and for simplicity we use  $w(n,m)=v(x_n,x_m)$.
Here $v$ excessive means that   $w(n,m)$ is weakly decreasing in $n$, and $v$ depressive means that $w(n,m)$ is weakly increasing in $m$. 
The meaning of $v$ balanced is the following: for each $n$ and $m$, $w(n,m)$ is the value of the matrix game (``local game" at $(n,m)$):
$$\left(\begin{array}{cc}
 w({n+1,m}) &  w({n+1,m+1})\\
 w({n,m}) &  w({n,m+1})\\
\end{array}
\right).$$ 
 
Clearly, if $v$ is excessive and depressive it is balanced, but proposition \ref{pro2} tells that {\it if $v$ is continuous}, then the converse also holds: balancedness implies excessiveness and depressiveness. The idea of the proof of proposition \ref{pro2} can be seen here as follows.

Suppose $v$ is balanced, but not excessive. Then one can find $n$ and $m$ such that $w({n+1},m)>w(n,m)$.
Because $w(n,m)$ is the value of the local game at $(n,m)$, we necessarily have $w(n+1,m+1)\leq w(n,m)$, and $w(n+1,m)>w(n+1,m+1)$. Consider now the ``local  game" at $(n+1,m)$. $w(n+1,m)$ is the value of the matrix:
 $$\left(\begin{array}{cc}
 w({n+2,m}) &  w({n+2,m+1})\\
 w({n+1,m}) &  w({n+1,m+1})\\
\end{array}
\right).$$ \noindent Since  $w(n+1,m)>w(n+1,m+1)$, we obtain  $w(n+2,m+1)\geq w(n+1,m).$
We have obtained: $w(n+2,m+1)\geq w(n+1,m)>w(n,m)\geq w(n+1,m+1)$, so
$$w(n+2,m+1)-w(n+1,m+1)>w(n+1,m)-w(n,m).$$
Iterating  the argument, we obtain that for each $p$,
$$w(n+p+1,m+p)\geq w(n+1,m)-w(n,m)>0.$$
And this is a contradiction with $w$ being continuous at infinity.\\

To conclude with this example, consider   the simple  case where the running payoff is given by $u(x,y)=|x-y|$. Player 1 wants to be far from Player 2, and Player 2 wants to be close to Player 1. If initially $n<m$, it is optimal for each player not to move, so $w(n,m)=|x_n-x_m|$. Suppose on the contrary that initially $n\geq m$, so that Player 1 is more to the right than Player 2. Then Player 2 has a simple optimal strategy which is to move  to the right if the current positions satisfy $x>y$, and to stay at $y$ if $x=y$. No matter how large is the initial difference $n-m$, Player 2 will succeed in being close to player, so that $w(n,m)=0$ if $n\geq m$.

 \section{Proof  of Theorem \ref{thm1},  part  1.}
 
 We prove here  propositions  \ref{pro2}, \ref{pro2.5} and  \ref{pro3}.

 \subsection{Proof of proposition \ref{pro2}}

By symmetry, suppose that $\Gamma$ is weakly acyclic and $\Lambda$ Leavable.  Let us prove that any balanced continuous function $v$ is excessive-depressive. First let us first prove that $v$ is excessive. 

 Fix any $(x_0,y_0)$ in $X\times Y$, and  $p_1 \in \Gamma(x_0)$. A direct consequence of balancedness is the existence of $q_1$ in $\Lambda(y_0)$ such that $v(p_1,q_1)\leq v(x_0,y_0)$. Now, $p_1$ is in $\Delta(X)$ and $q_1$ is in $\Delta(Y)$. One has to be careful that there may not exist $p_2\in  \tilde{\Gamma}(p_1)$ such that for all $q_2\in \tilde{\Lambda}(q_1)$, $v(p_2,q_2)\geq v(p_1,q_1)$. This is because, $\tilde{v}$ being affine in each variable,  $v(p_1,q_1)$ can be interpreted as   the value of the auxiliary game where first $x$ and $y$ are chosen according to $p_1\otimes q_1$ and {\it observed by the players}, then players respectively choose $p\in \Gamma(x)$ and $q\in \Lambda(y)$ and finally Player 1's payoff is $v(p,q)$. And to play well in this game Player 1 has to know the realization of $q_1$ before choosing $p$. However since $y_0$ is a Dirac measure,  balancedness implies  that there exists $p_2\in \tilde{\Gamma}(p_1)$ such that $v(p_2,q_1)\geq v(p_1,y_0).$ We have obtained the following lemma:
 
 \begin{lem}\label{lem2}
 Given $(x_0,y_0)$ in $X\times Y$, and  $p_1 \in \Gamma(x_0)$, there exists  $q_1$ in $\Lambda(y_0)$ and $p_2\in \tilde{\Gamma}(p_1)$ such that:
 $v(p_1,q_1)\leq v(x_0,y_0)\;{\rm and}\; v(p_2,q_1)\geq v(p_1,y_0).$
 \end{lem}

\vspace{0,5cm}

We now prove the proposition. Define, for   $x$ in $X$, $$h(x)=\Max\{v(p,y)-v(x,y), y\in Y, p \in \Gamma(x)\}.$$

\noindent  $h$ is continuous. We put $Z=\argmax_{x\in X} h(x)$, and consider $x_0\in \argmin_{x \in Z} \varphi(x)$, where $\varphi$ comes from the definition of $\Gamma$ weakly acyclic. 
 
$x_0 \in Z$, so $v$ is excessive  if and only if $h_0=_{def}h(x_0)\leq 0$.
By definition of $h(x_0)$, there exists $p_1\in \Gamma(x_0)$ and $y_0$ in $Y$ such that $v(p_1,y_0)-v(x_0,y_0)=h_0.$\\

By lemma \ref{lem2}, there exists $q_1$ in $\Lambda(y_0)$ and $p_2\in \tilde{\Gamma}(p_1)$ such that $v(p_1,q_1)\leq v(x_0,y_0)$ and $v(p_2,q_1)\geq v(p_1,y_0).$ Consequently, 
$$v(p_2,q_1)-v(p_1,q_1)\geq v(p_1,y_0)-v(x_0,y_0)=h_0.$$
One can now  find $y_1$ in $Y$ such that $v(p_2,y_1)-v(p_1,y_1)\geq h_0$, and since $p_2\in \tilde{\Gamma}(p_1)$ it implies that $\Supp(p_1) \subset Z$. But $p_1\in  \Gamma(x_0)$, so by  definition of $x_0$ and weak acyclicity, we obtain that $p_1=\delta_{x_0}.$ So $h_0=0$, and $v$ is excessive. \\

Let us now prove that when $v$ is balanced and excessive then it is depressive. For every $(x,y)$, and every $p\in \Gamma(x)$, we have  $v(x,y)\geq v(p,y)$. Thus, for every $(x,y)$,  $p\in \Gamma(x)$ and every $q \in\Lambda(y)$, $v(x,q)\geq v(p,q)$ and consequently,  $\min_{q \in\Lambda(y)}v(x,q)\geq \min_{q \in\Lambda(y)} v(p,q)$. Taking the maximum in $p\in \Gamma(x)$ and using that $v$ is balanced implies that $\min_{q \in\Lambda(y)}v(x,q)\geq v(x,y)$. Since $\Lambda$ is leavable, we have equality and so $v$ is depressive  with respect to $Y$.

\subsection{Proof of proposition \ref{pro2.5}} Let $(\lambda_n)_n$ be a vanishing sequence of discount factors such that $\|v_{\lambda_n}-v\|\rightarrow_{n \to \infty}0$. Fix $(x,y)$ in $X\times Y$, by symmetry it is enough to show  that there exists $ p \in \Gamma^\infty(x)$ such that  $v(x,y)\leq v(p,y)\leq u(p,y)$. If $v(x,y)\leq u(x,y)$, it is enough to consider $p=\delta_x$, so we assume  $v(x,y)> u(x,y)$. For $n$ large enough, $v_{\lambda_n}(x,y)>u(x,y)+ \lambda_n$.\\

Fix $n$. We   define inductively a sequence $(p_t^n)_{t=0,...,T_n}$ in $\Delta(X)$, with $T_n\geq 1$,  by: 

1) $p_0^n=\delta_x$,

2)   for each $t\geq 0$ such that  $v_{\lambda_n}(p_t^n,y)>u(p_t^n,y)+ \lambda_n$, we define $p_{t+1}^n$ in $\tilde{\Gamma}(p_t^n)$ by: $$p_{t+1}^n \in \argmax_{p \in \tilde{\Gamma}(p_t^n)} \left(\lambda_n u(p_t^n,y)+ (1-\lambda_n) v_{\lambda_n}(p,y)\right).$$
We have $\max_{p \in \tilde{\Gamma}(p_t^n)} \left(\lambda_n u(p_t^n,y)+ (1-\lambda_n) v_{\lambda_n}(p,y)\right)\geq v_{\lambda_n}(p_t^n,y)$, so:
\begin{equation}\label{eq1} \lambda_n u(p_t^n,y)+ (1-\lambda_n) v_{\lambda_n}(p_{t+1}^n,y)\geq v_{\lambda_n}(p_t^n,y).\end{equation}
Since $u(p_t^n,y)< v_{\lambda_n}(p_t^n,y) - \lambda_n$, we obtain:
\begin{equation} \label{eq2}v_{\lambda_n}(p_{t+1}^n,y)\geq v_{\lambda_n}(p_t^n,y)+ \frac{{\lambda^2_n}}{1-\lambda_n}>v_{\lambda_n}(p_t^n,y).\end{equation}
Since  $\frac{\lambda^2_n }{1-\lambda_n}>0$ and $v_{\lambda_n}$ is bounded, there exists a first integer  $t=T_n$ where $v_{\lambda_n}(p_{T_n}^n,y)\leq u(p_{T_n}^n,y)+ \lambda_n$, and we stop here the definition of the sequence   $(p_t^n)_{t=0,...,T_n}$. Inequalities (\ref{eq1}) and (\ref{eq2}) above give: \begin{equation} \label{eq3}
v_{\lambda_n}(x,y)\leq   v_{\lambda_n}(p_{T_n-1}^n,y)\leq \lambda_n u(p_{T_n-1}^n,y)+ (1-\lambda_n) v_{\lambda_n}(p_{T_n}^n,y).
\end{equation}

Define now $p^n=p^n_{T_n}$ for each $n$. $p^n\in \tilde{\Gamma}^{T_n}(x)\subset \Gamma^\infty(x)$ for each $n$, and we consider a limit point $p^*\in  \Gamma^\infty(x)$ of   $(p^n)_n$. Because $v_{\lambda_n}$ is an equicontinuous family   converging to $v$,   we obtain the convergence (along a subsequence) of $v_{\lambda_n}(p_{T_n}^n,y)$ to $v(p^*,y)$. 
 Passing to the limit in the inequality defining $T_n$ then gives:

$$v(p^*,y)\leq u(p^*,y).$$
Finally, passing to the limit in    (\ref{eq3}) shows: 
$v(x,y)\leq v(p^*,y).$

%$v_{\lambda_n}(p_{T_n}^n,y)\xrightarrow[n\to \infty]{} .

\subsection{Proof of proposition \ref{pro3}} 

We start with a lemma.

\begin{lem} \label{lem3}
Assume $v$ in $\mathcal{C}(X\times Y)$ is excessive. Then for all $(x_0,y_0)$ in $X\times Y$ and all $p\in \Gamma^{\infty}(x_0)$, we have: $v(p,y_0)\leq v(x_0,y_0).$
\end{lem}

\noindent{\bf Proof:} $p=\lim_n p_n$, with $p_n\in \tilde{\Gamma}^n(x_0)$ for each $n$. It is enough to prove that $v(p_n,y_0)\leq v(x_0,y_0)$ for each $n$, and we do the proof by induction on $n$. The case $n=1$ is clear by definition of $v$ excessive. Since $p_{n+1}\in \tilde{\Gamma}(p_n)$, it is enough to prove that for $p'$ in $\Delta(X)$ and $p''\in \tilde{\Gamma}(p')$, we have $v(p'',y_0)\leq v(p',y_0).$  By definition of $\tilde{\Gamma}$, $(p',p'')$ is in the closure of $\conv (Graph\Gamma)$. 

$y_0$ is fixed, and the function $h: p \longrightarrow v(p,y_0)$ is affine continous on $\Delta(X)$. The set $$D=_{def} \{(p',p'')\in \Delta(X)\times \Delta(X), h(p'')\leq h(p')\}.$$ is convex and compact, and we want to show that $Graph (\tilde{\Gamma})=\overline{\conv} (Graph\Gamma)\subset D$. It it enough to prove that $Graph(\Gamma)\subset D$, and this is implied by the fact that $v$ is excessive. This concludes the proof of lemma \ref{lem3}.\\

We now prove the proposition. Assume one of the gambling houses is strongly acyclic, and let $v_1$ and $v_2$ satisfying the conditions of proposition \ref{pro3} (are continuous, balanced, $v_1$ satisfies $P1$ and $v_2$ satisfies $P2$). We will show that $v_1\leq v_2$. 

By symmetry, suppose that $\Gamma$ is strongly  acyclic.  From Proposition \ref{pro2}, $v_1$ and $v_2$ are excessive (in $X$) and depressive (in $Y$). $v_1-v_2$ being continuous on   $X\times Y$, define the compact set:
$$Z={\rm Argmax}_{(x,y)\in X\times Y} \; v_1(x,y)-v_2(x,y).$$
\indent Consider now $\varphi$ u.s.c. given by the strong acyclicity   condition of $\Gamma$. The set $Z$ being compact, there exists $(x_0,y_0)$ minimizing  $\varphi(x)$ for $(x,y)$ in $Z$.

By $v_2$ satisfying $P2$, there exists $q$ in $\Lambda^{\infty}(y_0)$ such that $v_2(x_0,y_0)=v_2(x_0,q)\geq u(x_0,q)$. Thus, there is $y'_0 \in \Supp(q) $ such that  $v_2(x_0,y'_0)\geq u(x_0,y'_0)$.

Because $v_1$ is depressive, by lemma \ref{lem3} we have $v_1(x_0,q)\geq v_1(x_0,y_0)$ and we obtain $v_1(x_0,q)-v_2(x_0,q)\geq v_1(x_0,y_0)-v_2(x_0,y_0)$. Since $(x_0,y_0)$ is in $Z$, $\{x_0\}\times \Supp(q) \subset Z$. Thus, $(x_0,y'_0)\in Z$. Obviously,  $(x_0,y'_0)$ also minimizes  $\varphi(x)$ for $(x,y)$ in $Z$ (the minimum value remains unchanged: $\varphi(x_0)$).

By $v_1$ satisfying $P1$, there exists $p$ in $\Gamma^{\infty}(x_0)$ such that: $v_1(x_0,y'_0)=v_1(p,y'_0)\leq u(p,y'_0)$. Because $v_2$ is excessive, by lemma \ref{lem3} we have $v_2(p,y'_0)\leq v_2(x_0,y'_0)$ and we obtain $v_1(p,y'_0)-v_2(p,y'_0)\geq v_1(x_0, y'_0)-v_2(x_0,y'_0)$. Since $(x_0,y'_0)$ is in $Z$,  $\Supp(p)\times \{y'_0\} \subset Z$. By definition of $(x_0,y'_0)$, this implies:
$\varphi(p)\geq \varphi(x_0)$. The definition of $\varphi$ now gives that $p$ is the Dirac measure on $x_0$. We obtain that   $v_1(x_0,y'_0)=v_1(p,y'_0)\leq u(x_0,y'_0)$. So $\Max_{(x,y)\in X\times Y} v_1(x,y)-v_2(x,y)= v_1(x_0,y'_0)-v_2(x_0,y'_0)\leq u(x_0,y'_0)-u(x_0,y'_0)=0$, and thus $v_1\leq v_2$. \\

%\begin{rem} \rm
%The proof of the proposition shows that any excessive-depressive continuous function satisfying $P1$ is smaller than any excessive-depressive continuous function satisfying $P2$. 
%\end{rem}

\section{A weakly acyclic  game without limit value}\label{contrex}

We conclude here the proof of theorem \ref{thm1} by providing a counterexample to the convergence of $(v_\lambda)$ in a weakly acyclic non expansive gambling house.

\subsection{The counter-example}

The states and transitions for  Player 1 are  as in example \ref{exa2}:   \vspace{1cm}

    \begin{center}
\setlength{\unitlength}{1mm}
\begin{picture}(60,30)

 \put(30,30){\circle{12}}
  \put(29,32){a}

     \put(10,5){\circle{12}}
  \put(9,7){b}

         \put(50,5){\circle{12}}
  \put(49,7){c}

    \textcolor{blue} {      { \qbezier(25,25)(23,20)(12,8)}}
 \textcolor{blue} {   { \put(19,18){\vector(-1,-1){1}}}}
 
     \textcolor{blue} {      { \qbezier(30,30)(33,45)(23,35)}}
 \textcolor{blue} {   { \put(27,39){\vector(-1,0){1}}}}
 
 \textcolor{blue} {     \put(10,15){$\alpha_a$}}
  \textcolor{blue} {     \put(30,35){$1-\alpha_a- \alpha_a^2$}}
   
    \textcolor{blue} {     \put(29,12){$\alpha_a^2$}}
 
 \textcolor{blue} {     { \qbezier(22,26)(37,20)(40,9)}}
 \textcolor{blue} {   { \put(34,17){\vector(1,-1){1}}}}

     \textcolor{red} {      { \qbezier(-5,8)(10,30)(15,30)}}
 \textcolor{red} {   { \put(0,18){\vector(1,1){1}}}}
  \textcolor{red} {     \put(-5,20){$\alpha_b$}}

     \textcolor{red} {      { \qbezier(-10,0)(-15,-15)(-5,0)}}
 \textcolor{red} {   { \put(-10,-5){\vector(1,1){1}}}}

  \textcolor{red} {     \put(-20,-10){$1-\alpha_b- \alpha_b^2$}}

 \textcolor{red} {     { \qbezier(-10,5)(15,5)(25,5)}}
 \textcolor{red} {   { \put(15,5){\vector(1,0){1}}}}
     \textcolor{red} {     \put(10,-5){$\alpha_b^2$}}
\end{picture}
    
      \end{center}

\vspace{1cm}

\noindent The set of states of Player 1 is  $X=\{a,b,c\}$. The difference with  example \ref{exa2} is that the set  of possible choices for $\alpha_a$ and $\alpha_b$ may be smaller than $[0,1/2]$.  Here   $\alpha_a$ and $\alpha_b$ now belong to some fixed compact set $I\subset [0,1/2]$  such\footnote{Except in part 2 of theorem \ref{thm2} where we will consider the case $I=\{0,1/4\}$.}  that 0 is in the closure of $I\backslash\{0\}$.  
Then $0\in I$, the transitions are  leavable and non expansive.
  States $a$ and $b$ are  symmetric, this  gambling house  is weakly acyclic, with $\varphi(a)= \varphi(b)=1, \varphi(c)=0$, but not strongly acyclic  since    $a\in \Gamma^\infty(b)$ and $b\in \Gamma^\infty(a)$.

  The gambling house of Player 2 is a  copy of the gambling house  of Player 1, with state space  $Y=\{a',b',c'\}$ and a compact set of choices $J\subset [0,1/2]$   such that  0  is  in the closure of $J\backslash\{0\}$. The  unique difference between the gambling houses of the players is that 
 $I$ and $J$ may be different. Payoffs are   simple:
      $$u(x,y)=0\; {\rm  if}\;  x=y, \; \; u(x,y)=1 \; {\rm if}\;  x\neq y.$$
      The $u$ function can be written as follows     \; $
\begin {tabular}{ccccccc}
& \multicolumn{1}{c}{$$} &\multicolumn{1}{c}{$a'$} &\multicolumn{1}{c}{$b'$}  &\multicolumn{1}{c}{$c'$}  \\
\cline{3-5}
& \multicolumn{1}{c}{$a$} &\multicolumn{1}{|c|}{$0$} &
\multicolumn{1}{|c|}{$1$}  &
\multicolumn{1}{|c|}{$1$}\\
\cline{3-5}
& \multicolumn{1}{c}{$b$} &\multicolumn{1}{|c|}{$1$} &
\multicolumn{1}{|c|}{$0$}  &
\multicolumn{1}{|c|}{$1$}\\
\cline{3-5}
& \multicolumn{1}{c}{$c$} &\multicolumn{1}{|c}{$1$} &
\multicolumn{1}{|c|}{$1$}  &
\multicolumn{1}{|c|}{$0$}\\
\cline{3-5}
\end{tabular}
$\;, with a clear  interpretation : Player 1 and Player 2 both move on a  space with 3 points, Player 2 wants to be at the same location as Player 1, and Player 1 wants the opposite.\\

Here the gambling game is weakly acyclic but not strongly acyclic, and the following lemma shows that the uniqueness property of proposition \ref{pro3} fails.    

   \begin{lem}   \label{lem5}       Let   $v:X\times Y\longrightarrow \R$. Then A) and B) are equivalent: 
      
A) $ v \; {\it is} \;  {\it excessive},   \;  {\it depressive},$ and satisfies $P1$ and $P2$,
% $$\forall (x,y)\in X\times   Y, \exists p \in \Gamma^\infty(x), v(x,y)=v(p,y)\leq u(p,y),$$
%$$\forall (x,y)\in X\times   Y, \exists q \in \Lambda^\infty(y), v(x,y)=v(x,q)\geq u(x,q). $$

B) There exists   $x\in [0,1]$ such that $v$ can be written:
    $$ 
\begin {tabular}{ccccccc}
& \multicolumn{1}{c}{$$} &\multicolumn{1}{c}{$a'$} &\multicolumn{1}{c}{$b'$}  &\multicolumn{1}{c}{$c'$}  \\
\cline{3-5}
& \multicolumn{1}{c}{$a$} &\multicolumn{1}{|c|}{$x$} &
\multicolumn{1}{|c|}{$x$}  &
\multicolumn{1}{|c|}{$1$}\\
\cline{3-5}
& \multicolumn{1}{c}{$b$} &\multicolumn{1}{|c|}{$x$} &
\multicolumn{1}{|c|}{$x$}  &
\multicolumn{1}{|c|}{$1$}\\
\cline{3-5}
& \multicolumn{1}{c}{$c$} &\multicolumn{1}{|c}{$0$} &
\multicolumn{1}{|c|}{$0$}  &
\multicolumn{1}{|c|}{$0$}\\
\cline{3-5}
\end{tabular}
$$
\end{lem}

\noindent{\bf Proof:} Assume $v$ satisfies $A)$. Since $u$ takes values in $[0,1]$, so does also $v$. 
 Because $c$ and $c'$ are absorbing,  we have $v(c,c')=u(c,c')=0$.  
 
 Consider now $v(c,a')$. $v$ being depressive, for any fixed $\beta^*>0$ in $J$ we have:
 $$v(c,a')\leq \beta^* v(c,b')+ {\beta^*}^2 v(c,c') + (1-\beta^*-{\beta^*}^2) v(c,a').$$
 \noindent and we obtain $v(c,a')\leq \frac{1}{1+\beta^*} v (c,b')$. But symmetrically we also have $v(c,b')\leq \frac{1}{1+\beta^*}v(c,a')$, and we get $v(c,a')=v (c,b')=0$. 
 
Consider now $v(a,c')$. By $P2$, we obtain that $v(a,c')\geq u(a,c')=1$, and $v(a,c')=1$. Similarly, $v(b,c')=1$.

Consider now $v(a,a')$. $v$ being excessive, for any   $\alpha>0$ in $I$ we have:
 $$v(a,a')\geq \alpha  v(b,a')+ {\alpha}^2 v(c,a') + (1-\alpha-{\alpha}^2) v(a,a').$$
 \noindent Hence $v(a,a')\geq \frac{1}{1+\alpha} v(b,a')$, and by assumption on $I$   we obtain $v(a,a')\geq v(b,a')$. By symmetry of the transitions between $a$ and $b$, $v(a,a')=v(b,a')$. Similarly, $v(a,b')=v(b,b')$.
 
 It only remains  to prove that $v(a,a')=v(a,b')$. $v$ being depressive, for any $\beta>0$ in $J$,
 $v(a,a')\leq \beta v(a,b') + \beta^2 v(a,c') +(1-\beta-\beta^2) v(a,a').$ By assumption on $J$, we get $v(a,a')\leq v(a,b')$. By symmetry of the transitions between $b$ and $b'$, $v(a,b')=v(a,a')$, and $v$ satisfies  $B)$. 
 
 One can easily check that $B)$ implies $A)$, and 
 the proof of lemma \ref{lem5} is complete.\\

The second part of theorem \ref{thm1} is a direct consequence of  the following result.

\begin{thm} \label{thm2}  $\;$

 1) If   $I= J=[0,1/4]$, the limit value exists and is:
$$
\begin {tabular}{ccccccc}
& \multicolumn{1}{c}{$$} &\multicolumn{1}{c}{$a'$} &\multicolumn{1}{c}{$b'$}  &\multicolumn{1}{c}{$c'$}  \\
\cline{3-5}
& \multicolumn{1}{c}{$a$} &\multicolumn{1}{|c|}{$1/2$} &
\multicolumn{1}{|c|}{$1/2$}  &
\multicolumn{1}{|c|}{$1$}\\
\cline{3-5}
& \multicolumn{1}{c}{$b$} &\multicolumn{1}{|c|}{$1/2$} &
\multicolumn{1}{|c|}{$1/2$}  &
\multicolumn{1}{|c|}{$1$}\\
\cline{3-5}
& \multicolumn{1}{c}{$c$} &\multicolumn{1}{|c}{$0$} &
\multicolumn{1}{|c|}{$0$}  &
\multicolumn{1}{|c|}{$0$}\\
\cline{3-5}
\end{tabular}
$$

2)  If     $J=[0,1/4]$ and $I=\{0,1/4\}$, the limit value exists and is:
$$
\begin {tabular}{ccccccc}
& \multicolumn{1}{c}{$$} &\multicolumn{1}{c}{$a'$} &\multicolumn{1}{c}{$b'$}  &\multicolumn{1}{c}{$c'$}  \\
\cline{3-5}
& \multicolumn{1}{c}{$a$} &\multicolumn{1}{|c|}{$0$} &
\multicolumn{1}{|c|}{$0$}  &
\multicolumn{1}{|c|}{$1$}\\
\cline{3-5}
& \multicolumn{1}{c}{$b$} &\multicolumn{1}{|c|}{$0$} &
\multicolumn{1}{|c|}{$0$}  &
\multicolumn{1}{|c|}{$1$}\\
\cline{3-5}
& \multicolumn{1}{c}{$c$} &\multicolumn{1}{|c}{$0$} &
\multicolumn{1}{|c|}{$0$}  &
\multicolumn{1}{|c|}{$0$}\\
\cline{3-5}
\end{tabular}
$$

3) If $J=[0,1/4]$ and $I=\{\frac{1}{4^{n}}, n \in \N^*\} \cup \{0\}$, then $v_\lambda$ diverges.  \end{thm}

\subsection{Proof of theorem \ref{thm2}}

We start with considerations  valid for the 3 cases of the theorem. We fix $J=[0,1/4]$ in all the proof, and only assume for the moment that $I$ is a compact subset of $[0,1/4]$ containing $0$ and $1/4$.

 Consider $\lambda\in(0,1)$. It is clear that $v_\lambda(c,c')=0$, and $v_\lambda(a,c')=v_\lambda(b,c')=1$. By symmetry of the payoffs and transitions, we have $v_\lambda(a,a')=v_\lambda(b,b')$,   $v_\lambda(b,a')=v_\lambda(a,b')$ and $v_\lambda(c,c')=v_\lambda(c,b')$, so we can write $v_\lambda$ as:
 $$ 
\begin {tabular}{ccccccc}
& \multicolumn{1}{c}{$$} &\multicolumn{1}{c}{$a'$} &\multicolumn{1}{c}{$b'$}  &\multicolumn{1}{c}{$c'$}  \\
\cline{3-5}
& \multicolumn{1}{c}{$a$} &\multicolumn{1}{|c|}{$x_\lambda$} &
\multicolumn{1}{|c|}{$y_\lambda$}  &
\multicolumn{1}{|c|}{$1$}\\
\cline{3-5}
& \multicolumn{1}{c}{$b$} &\multicolumn{1}{|c|}{$y_\lambda$} &
\multicolumn{1}{|c|}{$x_\lambda$}  &
\multicolumn{1}{|c|}{$1$}\\
\cline{3-5}
& \multicolumn{1}{c}{$c$} &\multicolumn{1}{|c}{$z_\lambda$} &
\multicolumn{1}{|c|}{$z_\lambda$}  &
\multicolumn{1}{|c|}{$0$}\\
\cline{3-5}
\end{tabular}
$$
with $x_\lambda$, $y_\lambda$ and $z_\lambda$ in $(0,1)$.

 $z_\lambda$ is indeed easy to compute.  If the game is at $(c,a')$, Player 1 can not move, and Player 2 wants to reach  $c'$ as fast as possible, so he will choose $\beta=1/4$    and we have (see definition \ref{def1}):
$z_\lambda=\lambda 1 +(1-\lambda) (\frac{1}{16} 0 +\frac{15}{16} z_\lambda),$
 so that: \begin{equation}  \label{eq0}  z_\lambda=\frac{ 16 \lambda}{1+15\lambda}\leq 16 \lambda.\end{equation}

\begin{pro} \label{pro11}  Assume $J=[0,1/4]$, $\min I=0$ and $\max I=1/4$. Then for $\lambda$ small enough,   \begin{equation} \label{eq31}z_{\lambda}<x_{\lambda}<y_{\lambda},\end{equation}
  \begin{equation} \label{eq32}
   \lambda \; x_\lambda = (1-\lambda) \max_{ \alpha \in I} \left(\alpha(y_\lambda-x_\lambda) + \alpha^2(z_\lambda-x_\lambda)\right),\end{equation}
  \begin{equation} \label{eq33}
 \lambda \; y_\lambda = \lambda +(1-\lambda)  \min_{\beta \in J}\left(\beta (x_\lambda -y_\lambda) +\beta^2 (1-y_\lambda)\right).\end{equation}

 \end{pro}

(\ref{eq32}) express the fact that at $(a,a')$ or $(b,b')$, it is optimal for Player 2 to play the pure strategy $\beta=0$ (stay at the same location and   wait until Player 1 has moved), and Player 1 can play a pure strategy $\alpha$ there. Similarly,  (\ref{eq33}) express the fact that at $(a,b')$ or $(b,a')$, it is optimal for Player 1 not to move. In spite of these simple intuitions,  the proof of the proposition is rather technical, and relegated to the appendix. \\

Taking for  granted proposition \ref{pro11}, we now proceed to the proof of theorem \ref{thm2}.  It is   simple to study the simple maximization problem  of Player 2 given by (\ref{eq33}), which is simply minimizing a concave polynomial on the interval $J=[0,1/4]$. 

  If $y_\lambda-x_\lambda\leq \frac{1}{2}(1-y_\lambda)$, the minimum  in  (\ref{eq33}) is achieved for $\beta=\frac{y_\lambda-x_{\lambda}}{2(1-y_\lambda)}$, and otherwise it is achieved for $\beta=\frac{1}{4}$. Hence for $\lambda$ small enough:
  
  \begin{equation}\label{eq14}
 {\rm If }\; y_\lambda-x_\lambda\leq \frac{1}{2}(1-y_\lambda), \;\; 4 \lambda(1-y_\lambda)^2=(1-\lambda)(y_\lambda-x_\lambda)^2.
 \end{equation}
    \begin{equation}\label{eq15}
 {\rm If }\; y_\lambda-x_\lambda\geq  \frac{1}{2}(1-y_\lambda), \;\;  (1-y_\lambda)(1+15\lambda)=4(1-\lambda) (y_\lambda-x_\lambda).
 \end{equation}
\noindent Notice that the 2 inequalities of   (\ref{eq15})  imply $1+15 \lambda \geq 2(1-\lambda) $, which is is not possible for $\lambda$ small.  Consequently, for   small discount factors the right hand side of (\ref{eq14}) holds, and we have proved, for $\lambda$ small enough,  the main equality of the proof:
 \begin{equation}\label{eq36}
   4 \lambda(1-y_\lambda)^2=(1-\lambda)(y_\lambda-x_\lambda)^2.   \end{equation}

This clearly implies:
$$y_\lambda-x_\lambda\longrightarrow_{\lambda \to 0}0.$$   For each $\lambda>0$, denote by $\alpha_\lambda\in I$ a maximizer in the   expression  (\ref{eq32}), so that 
\begin{equation}\label{eq37}\lambda x_\lambda=(1-\lambda)   \alpha_\lambda (y_\lambda-x_\lambda)+(1-\lambda) \alpha_\lambda^2(z_\lambda-x_\lambda). \end{equation} 
And since $x_\lambda>0$, $\alpha_\lambda>0$.\\
 
 The following lemma implies  part 2) of theorem \ref{thm2}.
 
  \begin{lem} \label{lem9} If 0 is an isolated point in $I$, then    $y_{\lambda}$ and $x_{\lambda}$ converge to 0.
\end{lem}

\noindent{\bf Proof:} In this case there exists $\alpha^*>0$ such that $\alpha_\lambda\geq \alpha^*$ for all $\lambda$. Passing to the limit in (\ref{eq37}) gives the result.\\

We will now prove parts 1) and 3) of the theorem.  The fact that  $I\subset J$ gives an advantage to Player 2, which can be quantified as follows. 

\begin{lem}\label{lem6} Assume that $y_{\lambda_{n}}$ and $x_{\lambda_{n}}$ converge to $v$ in $[0,1]$. Then $v\leq 1/2,$ and $y_{\lambda_n}-x_{\lambda_n} \sim 2 \sqrt{\lambda_n} (1-v).$
\end{lem}

\noindent{\bf Proof:} $(1-\lambda)(y_\lambda-x_\lambda)=2 \sqrt{\lambda}\sqrt{1-\lambda} (1-y_\lambda)$, so
$$x_\lambda(\lambda+\alpha_\lambda^2)= \lambda x_\lambda \alpha_\lambda^2+(1-\lambda) \alpha_\lambda^2 z_\lambda+2 \alpha_\lambda  \sqrt{\lambda} \sqrt{1-\lambda}(1-y_\lambda) \geq 2 \alpha_\lambda \sqrt{\lambda} x_\lambda,$$
since $\lambda+\alpha_\lambda^2- 2 \alpha_\lambda \sqrt{\lambda}\geq 0$. Dividing par $\alpha_\lambda \sqrt{\lambda}$ and passing to the limit gives $2(1-v)\geq 2 v$, so $v\leq 1/2,$ and the lemma is proved.\\

Consider again the concave optimization problem of Player 1 given  by equation (\ref{eq32}), and denote by $\alpha^*(\lambda)=\frac{y_\lambda-x_\lambda}{2(x_\lambda-z_\lambda)}>0$ the argmax of the unconstrained problem if Player 1 could choose any $\alpha\geq 0$. If $y_\lambda$ and $x_\lambda$ converge to $v>0$, then $\alpha^*(\lambda)\sim \sqrt{\lambda}\frac{1-v}{v}$, and Player 1 would like to play in the $\lambda$-discounted game at $(a,a')$  some $\alpha$ close to $\sqrt{\lambda}\frac{1-v}{v}$.

\begin{lem}  \label{lem7} Let $\lambda_n$ be a vanishing sequence of discount factors such that $\sqrt{\lambda_n}\in  I$ for each $n$. Then  $y_{\lambda_n}$ and $x_{\lambda_n}$ converge to 1/2. \end{lem}

\noindent{\bf Proof:} By considering a converging subsequence we can assume that $y_{\lambda_n}$ and $x_{\lambda_n}$ converge to some $v$ in $[0,1]$. By the previous lemma, $v\leq 1/2$, and we have to show that $v\geq 1/2$. We have for each $\lambda$ in the subsequence, since Player 1 can choose to play $\alpha=\sqrt{\lambda}$:
$$\lambda x_\lambda\geq (1-\lambda) \sqrt{\lambda}(y_\lambda-x_\lambda) +(1-\lambda) \lambda (z_\lambda-x_\lambda),$$
so $$x_\lambda(2-\lambda)\geq (1-\lambda)z_\lambda+ (1-\lambda)\frac{y_\lambda-x_\lambda}{\sqrt{\lambda}}.$$
By passing to the limit, we get $2v\geq 2(1-v)$, and $v\geq 1/2$.

\begin{lem} \label{lem8} Let $\lambda_n$ be a vanishing sequence of discount factors such that for each $n$,  the open interval $(\frac{1}{2}\sqrt{\lambda_n}, 2 \sqrt{\lambda_n})$ does not intersect $I$. Then    $\limsup_n y_{\lambda_n}\leq 4/9$. \end{lem}

\noindent{\bf Proof:} Suppose that (up to a subsequence) $x_{\lambda_n}$ and $y_{\lambda_n}$ converges to some $v\geq 4/9$. It is enough to show that $v=4/9$. We know that $v\leq 1/2$ by lemma \ref{lem6}, and since $\alpha^*(\lambda)\sim \sqrt{\lambda}\frac{1-v}{v}$ we have  $ \frac{1}{2}\sqrt{\lambda}\leq \alpha^*(\lambda)\geq 2\sqrt{\lambda}$ for $\lambda$ small in the sequence. By assumption $(\frac{1}{2}\sqrt{\lambda}, 2 \sqrt{\lambda})$ contains no point in $I$ and the objective function of Player 1 is increasing from 0    to $\alpha^*(\lambda)$ and decreasing after $\alpha^*(\lambda)$. There are 2 possible cases:

If  $\alpha_\lambda\leq \frac{1}{2}\sqrt{\lambda}$ we have: 
 $$\lambda x_\lambda \leq \frac{1}{2} (1-\lambda) \sqrt{\lambda}(y_\lambda-x_\lambda) + \frac{1}{4}(1-\lambda)\lambda (z_\lambda-x_\lambda).$$
Dividing by $\lambda$ and passing to the limit gives: $v\leq 1-v-\frac{1}{4}v$, i.e. $v\leq \frac{4}{9}.$

Otherwise, $\alpha_\lambda> 2 \sqrt{\lambda}$ and we have: 
 $$\lambda x_\lambda \leq 2 (1-\lambda) \sqrt{\lambda}(y_\lambda-x_\lambda) +4(1-\lambda)\lambda (z_\lambda-x_\lambda).$$
Again, dividing by $\lambda$ and passing to the limit gives: $v\leq 4(1-v)- {4}v$, i.e. $v\leq \frac{4}{9}.$\\

Finally, lemma \ref{lem7} proves  part  1) of theorem \ref{thm2}, whereas lemmas  \ref{lem7} and  \ref{lem8} together   imply part 3), concluding the proof of theorem \ref{thm2}.

\begin{rems} $\;$\rm

$\bullet$  In case 3) of theorem \ref{thm2}, it is  possible to show,  using   lemma  (\ref{lem8}), that $\liminf x_\lambda=\liminf y_\lambda=4/9$,  and $\limsup x_\lambda=\limsup y_\lambda=1/2$.

$\bullet$ It is not difficult to adapt lemma  (\ref{lem8})  to show the divergence of $(v_\lambda)$ as soon as $J=[0,1/4]$ and  $I$ satisfies: 

 a)  there exists a sequence  $(\lambda_n)$ converging to 0 such that   $\sqrt{\lambda_n}\in I$ for each $n$, and 
 
   b) there exist $\eta>0$ and a sequence  $(\lambda_n)$ converging to   0   such that for each $n$,  $I$ does not intersect the interval $ [\sqrt{\lambda_n}(1-\eta), \sqrt{\lambda_n}(1+\eta)]$.

$\bullet$ It is   important for the counterexample that $I=\{\frac{1}{2^{2n}}, n \in \N^*\} \cup \{0\}$ is not semi-algebraic. Indeed, it has been showed that if we assume  $X$ and $Y$   finite, and the transitions  $\Gamma$, $\Lambda$ and the payoff $u$  to be definable in some o-minimal structure, then $(v_\lambda)_\lambda$ converges  \cite{BGV_2015}.\end{rems}

  \section{Corollaries and Extensions}

%Recall that all gambling houses are supposed non-expensive and Leveable.

%\subsection{Other Evaluations}
%
%Since we have a stochastic game with continuous transitions with compact state space and where $v_{\lambda}$ uniformly converges, from a recent Tauberian theorem of Zilioto (2015), $v_n$ (the value of the $n$-finitely repeated game) also uniformly converges, to the same limit. 
%
%In fact, existence of an asymptotic value may be extended to more general classes of evaluations. More precisely, let $\theta = (\theta_t)$ be a general evaluation, that is probability distribution on $\mathbf{N}^*$ where $\theta_t$ is the weigh of stage $t$. Consider the gambling game where the total payoff of the game is given by $\sum_t \theta_t u(x_t,y_t)$. By Sion's minmax theorem, this game has a value denoted $v_{\theta}$. 
%
%The $n$-finitely repeated game corresponds to $\theta_t = \frac{1}{n}$ for $t=1,...n$ and $\theta_t = 0$, for $t>n$. The $\lambda$-discounted game corresponds to $\theta_t =\lambda (1-\lambda)^{t-1}$.
%
%An evaluation $\theta = (\theta_t)$ is called decreasing if for all $t\geq 1$, $\theta_t \geq \theta_{t+1}$.
%
%\begin{thm} \label{thm1-bis}
%Under the assumptions of theorem \ref{thm1}, any sequence $\theta^n$ of decreasing evaluations uniformly converges, as $\theta^n_1 \rightarrow 0$, to  $\lim_{\lambda \rightarrow 0} v_{\lambda}$, the unique excessive-depressive continuous functions satisfying P1 and P2.
%\end{thm}

\subsection{Gambling houses (or Markov Decision Processes)}

 We assume here that there is a unique player, i.e. that $Y$ is a singleton. Then non expansiveness is enough to guarantee the  uniform convergence of $(v_\lambda)_{\lambda}$ (as well as the uniform value, see    \cite{Renault_2011}) and the limit $v$ can be characterized  as follows \cite{RV_2013}: for all $x$ in $X$, 
 \begin{eqnarray*} v(x) &= &  \inf \large\{ 
w(x),  w:\Delta(X)\rightarrow \R \;  {\rm affine}\;  C^0  \; s.t.\\
 &\; &  (1)\;  \forall x'\in X, w(x')\geq \sup_{p \in \Gamma(x')} w(p)}, \; and\;   { (2) \;  \forall r \in R, w(r)\geq u(r)   \},
 \end{eqnarray*}
where  ${R=\{p \in \Delta(X), (p,p)\in {Graph \; \tilde{\Gamma}}\}}$  is interpreted  as   the set of {\it invariant measures} for  the gambling house (which is not necessarily  leavable here). If we  moreover  assume that the gambling house is leavable, then $R=\Delta(X)$ and we recover the  fundamental theorem  of gambling \cite{DS_1965, MS_1996}, namely, $(v_\lambda)$  uniformly converges  to:
$$v =\min\{w\in {\cal C}(X), w \;{\rm excessive}\; , w\geq u \} =\min\{w\in {\cal  B}(X), w\;{\rm excessive}\; , w\geq u \}.$$
\noindent It is also  easy to see that $v(x) = \sup_{p \in \Gamma^{\infty}(x)} u(p)$ for each $x$. \\

Our approach will lead to other  characterizations. We don't   assume any acyclicity condition in the following result.

\begin{thm} \label{thm6}
Consider a one player (leavable and non-expansive) gambling house. Then $(v_\lambda)$ uniformly converges to the unique function $v$ in $\mathcal{C}(X)$ satisfying: 
$v \; {\it is} \;  {\it excessive}, $
 $P1: \forall x\in X, \exists p \in \Gamma^\infty(x), v(x) = \tilde{v}(p)\leq \tilde{u}(p),$ and 
$P2: v\geq u.$
\end{thm}

\noindent \textbf{Proof:} From proposition \ref{pro2.5} any accumulation point of $(v_\lambda)$ is excessive and satisfies $P1$ and $P2$. Thus, we just need to show uniqueness, which is a direct consequence of the following lemma.

\begin{lem}
If $v_1 \in \mathcal{B}(X)$ satisfies $P1$ and   $v_2 \in \mathcal{C}(X)$ is excessive and satisfies $P2$, then $v_1 \leq v_2$. 
\end{lem}

\noindent \textbf{Proof:} Take any $x \in X$. Then there is $p \in \Gamma^\infty(x)$ such that $v_1(x)= v_1(p)\leq u(p)$. Because $v_2$ is excessive and continuous, by lemma \ref{lem3} $v_2(x) \geq v_2(p)$. Since  $v_2\geq u$, we have $v_2(p) \geq u(p)$. Consequently, one has $v_2(x)\geq u(p)$. Thus, $v_1(x)-v_2(x) \leq  0$, as desired.\\

 Using the gambling fundamental theorem, we obtain new viewpoints on the characterization of the limit value in leavable gambling houses. 

\begin{cor}\label{corthm6}
Consider a one player (leavable and non-expansive) gambling house. Then the asymptotic value exists and is:\\ 
(1) the smallest excessive function $v$ in $\mathcal{B}(X)$ satisfying $P2$;\\
(2) the largest excessive function $v$ in $\mathcal{B}(X)$ satisfying $P1$;\\
(3) the unique excessive function $v$ in $\mathcal{B}(X)$ satisfying $P1$ and $P2$.\\
Moreover, $v$ is continuous.
\end{cor}

%When we the house is strongly acyclic, we can prove other characterizations. 

%\textbf{Remark: on peut facilement generaliser cette characterization pour 1 joueur sans leavable ... mais laissons ca pour l'article suivant, quand on aura compris deux joueurs!}

\subsection{Equivalent characterizations}

\begin{defi} Given $g$ in ${\cal B}(X\times Y)$,
$Exc_{\Gamma}(g)$ is the smallest excessive (w.r.t. $X$) function not lower than $g$, and 
$Dep_{\Lambda}(g)$ is the largest   depressive (w.r.t. $Y$) function not greater than $g$.
\end{defi}

$Exc_{\Gamma}$ is usually called the {\it r\'eduite} operator and $Dep_{\Lambda}(g)=-Exc_{\Gamma}(-g)$. In splitting games, $Exc_{\Gamma}=Cav_{X}$  is the concavification operator on $X$ and $Dep_{\Lambda}(g)=Vex_{Y}$ is the convexification operator on $Y$. We introduce the following definition by analogy   with the Mertens-Zamir characterization.

\begin{defi} A function $v$ in ${\cal B}(X\times Y)$ 
satisfies the $MZ$-characterization if:
\begin{eqnarray*} MZ1: & \forall (x,y)\in X\times   Y, v(x,y)= Exc_{\Gamma} \min (u,v)(x,y), \\
 {\rm and}\; MZ2: & \forall (x,y)\in X\times   Y, v(x,y)= Dep_{\Lambda} \max (u,v)(x,y).
 \end{eqnarray*}
\end{defi}

We now introduce other properties, by analogy with the one established for splitting games and repeated games with incomplete information, see for instance \cite{Laraki_2001a, Laraki_2001b, RS_2001}. 

\begin{defi}  Let $v$ be  in ${\cal B}(X\times Y)$. 

1) For each $(x,y)$  in $X\times   Y$, 
\begin{center}
$x$ is extreme for $v(\cdot,y)$ if $\arg \max_{p \in \Gamma^{\infty}(x)} \widetilde{v}(p,y)=\{\delta_x\}$.\\
$y$ is extreme for $v(x,\cdot)$ if $\arg \min_{q \in \Lambda^{\infty}(y)} \widetilde{v}(x,q)=\{\delta_y\}$.
\end{center}

2) $v$  satisfies  the $E$-characterization if:
\begin{eqnarray*} E1: &{for\; all\;} (x,y)\in X\times   Y, {\it if \ x \ is \ extreme \ for \ v(\cdot,y) \ then \ } v(x,y) \leq u(x,y),\\
 {\rm and}\;   E2: &{for\; all\;} (x,y)\in X\times   Y, {\it if \ y \ is \ extreme \ for \ v(x,\cdot) \ then \ } v(x,y)\geq u(x,y).
 \end{eqnarray*}

\end{defi}

\begin{pro}
Consider a standard gambling game and let $v$ in $\mathcal{C}(X \times Y)$ be   excessive and depressive. Then:
\begin{center}
\hspace{0,7cm}
$v$ satisfies $MZ1$ $\implies$  $v$ satisfies  $P1$ $\implies$  $v$ satisfies  $E1$, \\
 
and $v$ satisfies $MZ2$ $\implies$  $v$ satisfies  $P2$ $\implies$  $v$ satisfies  $E2$. \\
\end{center}
\end{pro}

%This is not necessarily the case in theorem \ref{thm5} and corollary \ref{corthm5}.

\noindent \textbf{Proof.} Let $v$ be a continuous  excessive function that satisfies $MZ1$. Fix $y$ and define  for each $x$, $f(x)=\min (v(x,y),u(x,y))$. Then, for every $x$, $v(x,y)= Exc_{\Gamma} (f)(x)$. We consider the gambling house for Player 1 where the state of Player 2 is fixed to $y$ and the payoff is given by $f$. From corollary \ref{corthm6}, there is $p\in \Gamma^{\infty} (x)$ such that $v(x,y)=v(p,y)\leq f(p)$. Since $f(p) \leq u(p,y)$, $v$ satisfies $P1$. 

Now, let $v$ be an excessive continuous  function that satisfies $P1$. Take any $x$ and $y$ and suppose that $x$ is extreme for $v(\cdot,y)$. By $P1$, there is $p^*\in \Gamma^{\infty} (x)$ such that $v(x,y)=v(p^*,y)\leq u(p^*,y)$. Because $v$ is excessive and continuous, by lemma \ref{lem3} we have $p^* \in \arg \max_{p \in \Gamma^{\infty}} \widetilde{v}(p,y)$. Because $x$ is extreme for $v(\cdot,y)$, $p^*=\delta_x$ and so, $v(x,y) \leq u(x,y)$. Consequently, $E1$ is satisfied. 

By symmetry, $MZ2$ $\implies$   $P2$ $\implies$   $E2$.\\

\begin{rem} \rm It is easy to find examples where $E1$ is satisfied but $MZ1$ is not. For instance, assume that $Y$ is a singleton, and that $\Gamma(x)=\Delta(X)$ for each $x$ in $X$. Consider the constant, hence excessive, functions $u=0$ and $v=1$. $v$ has no extreme points hence satisfies $E1$, but $Exc_{\Gamma} \min (u,v)=u$ and $v$ does not satisfy $MZ1$. \end{rem}

\begin{pro} \label{thm7}
Consider a standard gambling game and let $v$ be an excessive-depressive function in $\mathcal{C}(X \times Y)$. Then:\\
($\Gamma$ strongly acyclic) and ($v$ satisfies $E1$) $\implies$ ($v$ satisfies $MZ1$),\\ 
and;\\
($\Lambda$ strongly acyclic) and ($v$ satisfies $E2$) $\implies$ ($v$ satisfies $MZ2$).

Consequently, if the gambling game is strongly acyclic, characterizations $MZ$, $P$ and $E$ are equivalents.
\end{pro}

\noindent \textbf{Proof.}  Let $v$ be excessive-depressive that satisfies $E1$. Fix $y \in Y$. We want to show that $v(x,y)=g(x,y)$ where $g=Exc_{\Gamma} (f)$ and $f =\min(u,v)$.  $g$ is continuous by corollary \ref{corthm6}. Since $v$ is excessive and $v\geq f$, we have $v\geq g$. Let $Z=\arg \max_{x\in X} v(x,y)-g(x,y)$ and let $x_0=\arg \min_{x \in Z} \varphi(x)$, where $\varphi$ comes form the definition of acyclicity. It is enough to prove that $v(x_0,y)\leq g(x_0,y)$. 

Suppose not. We have  $g(x_0,y) \geq f(x_0,y) $, so   (1) $ g(x_0,y) \geq u(x_0,y)$. Now, let $p_0\in \Gamma^{\infty}(x_0)$ such that $v(x_0,y)=v(p_0,y)$. Because $g$ is excessive and continuous, $g(p_0,y) \leq g(x_0,y)$ (lemma \ref{lem3}). Consequently, $v(p_0,y)-g(p_0,y) \geq v(x_0,y)-g(x_0,y)$. Consequently, $p_0$ is supported on $Z$. Thus, $\varphi(x_0,y) \leq \varphi (p_0,y)$, and by strong acyclicity  $p_{0}=\delta_{x_0}$. Thus $x_0$ is an extreme point of $v(\cdot,y)$. By $E_1$, we have (2) $v(x_0,y)\leq u(x_0,y)$.  By (1) and (2), $v(x_0,y) - g(x_0,y) \leq 0$.  A contradiction. \\

%Proof is in the appendix. In splitting games the $MZ$-characterization is the usual Mertens-Zamir system, $E$-characterization was proved in Laraki (2001), the $P$-characterization is new.  
%For idempotent games, the $P$-characterization takes a simplest form.

\subsection{Uniform optimal  strategies in   idempotent games} \label{subuniform}

\subsubsection{Idemptotent gambling games}
\begin{defi}
The gambling game is idempotent if $\Gamma \circ \Gamma= \Gamma$ and $\Lambda \circ \Lambda= \Lambda$.
\end{defi}

In that case, clearly $\Gamma=\Gamma^{\infty}$ and $\Lambda=\Lambda^{\infty}$. Any   state that could   be reached in several   stages can  be reached immediately in a single stage. This holds true for instance in splitting games. Notice also that   for any $\Gamma$, the multifunction $\Gamma^{\infty}$ is idempotent. If   $\Gamma \circ \Gamma= \Gamma$ , then  $\Gamma(x)=\tilde{\Gamma}^n(\delta_x)=\Gamma^\infty(x)$ for all $n$ and $x$  so if the gambling game is idempotent the notions of weak and strong acyclicity coincide.

An immediate corollary of theorem \ref{thm1} is the following. 

\begin{cor} \label{thm4}
Consider a standard   idempotent  gambling game where a  player has an acyclic gambling house.  Then $\{v_\lambda\}$ converges uniformly to the unique function $v$ in $\mathcal{C}(X \times Y)$ which is excessive, depressive and satisfies:
 $$Q1: \forall (x,y)\in X\times   Y, \exists p \in \Gamma(x), v(x,y) = v(p,y)\leq u(p,y),$$
$$Q2: \forall (x,y)\in X\times   Y, \exists q \in \Lambda(y), v(x,y) = v(x,q)\geq u(x,q).$$
\end{cor}

\subsubsection{Uniform value and optimal strategies}

In repeated and stochastic games, a stronger notion of limit value is given by the uniform value. As usual, a strategy of   Player 1, resp. Player 2, is a measurable rule giving at every stage $t$, as a function of past and current states, an element in  $\Gamma(x_t)$, resp. of  $\Lambda(y_t)$,  where $x_t$ and $y_t$ are the states of stage $t$. A pair $(x_1,y_1)$  of initial states and a pair of strategies $(\sigma, \tau)$ naturally define a probability on the set of plays $(X\times Y)^\infty$ (with the product $\sigma$-algebra, $X$ and $Y$ being endowed with their Borel $\sigma$-algebra), which expectation is written $\E_{(x_1,y_1),\sigma, \tau}$. \begin{defi} $w\in \mathcal{B}(X\times Y)$ is the uniform value of the gambling game and both players have optimal uniform strategies if: 
 
There exists  a strategy $\sigma$ of Player 1    that uniformly guarantees $w$: for any $\varepsilon >0$, there is $N$ such that for any  any $n\geq N$ and initial states $(x_1,y_1)$,   for any strategy $\tau$   of Player 2, $\E_{(x_1,y_1),\sigma, \tau} \left( \frac{1}{n} \sum_{t=1}^n u(x_t,y_t)\right)\geq w(x_1,y_1)-\varepsilon$. 

And similarly, there exists  a strategy $\tau$ of Player 2 that uniformly guarantees $w$: for any $\varepsilon >0$, there is $N$ such that for any $n\geq N$ and initial states $(x_1,y_1)$,   for any strategy   $\sigma$ of Player 1, $\E_{(x_1,y_1),\sigma, \tau} \left( \frac{1}{n} \sum_{t=1}^n u(x_t,y_t)\right)\leq w(x_1,y_1)+\varepsilon$. 
\end{defi}
\noindent It is known \cite{Sorin_2002} that the above conditions imply similar inequalities for discounted payoffs: for Player 1, the same strategy  $\sigma$ is such that for any $\varepsilon >0$, there is $\lambda_0>0$ satisfying:  for any $\lambda\leq \lambda_0$,  $(x_1,y_1)$ and $\tau$, $\E_{(x_1,y_1),\sigma, \tau} \left( \lambda \sum_{t=1}^\infty (1-\lambda)^{t-1} u(x_t,y_t)\right)\geq w(x_1,y_1)-\varepsilon$. And if the uniform value exists, it has to be $v=\lim_\lambda v_\lambda=\lim_n v_n$.

\subsubsection{Adapted strategies}

Our main theorem  \ref{thm1} suggests particularly interesting strategies for the players.
Consider again conditions $P1$ and $P2$, and fix a pair of states $(x,y)$. If $u(x,y) \geq v(x,y)$, the running payoff of Player 1 is at least as good as the payoff he should expect in the long run, so we may consider that Player 1 is ``quite happy" with the current situation and in order to satisfy $P1$ it is enough for him not to move, i.e. to choose $p=\delta_x$. If on the contrary $u(x,y) < v(x,y)$, Player 2 is   happy with the current situation and can choose $q=\delta_y$ to satisfy $P2$, whereas Player 1 should do something, and a possibility  is  to move towards a $p$ satisfying $P1$. This looks interesting for Player 1 because if Player 2 does not react, eventually the distribution on states will approach $(p,y)$ and (in expectation) Player 1 will be happy again with the current situation  since $u(p,y)\geq v(p,y)$.

\begin{defi}Let $w$ in ${\cal B}(X\times Y)$.

 A strategy  of Player 1 is adapted to $w$  if  whenever the current state is $(x,y)$, it  plays $p  \in \Gamma(x)$ such that $w(x,y)\leq w(p,y)\leq u(p,y)$.

  A strategy  of Player 2  is adapted to $w$  if  whenever the current state is $(x,y)$, it  plays $q \in \Gamma(y)$ such that $w(x,y)\geq w(x,q)\geq u(x,q)$.
\end{defi}

If $w$ satisfies $Q1$, resp. $Q2$,  Player 1, resp. Player 2,  has a strategy adapted to $w$ (using a measurable selection theorem \cite{AB_2006}). If moreover $w$ is excessive, we have $w(x,y)= w(p,y)\leq u(p,y)$. Mertens Zamir \cite{MertensZamir_1971} in repeated games with incomplete information and Oliu-Barton \cite{Miquel} in splitting games used similar strategies derived from the $MZ$-characterization instead of the $Q$-characterization. 

%In fact he proves that Player 1 (resp. Player 2) can uniformly guarantee any concave-convex function satisfying MZ1 (resp.

\subsubsection{Vanishing $L_1$-variation}

In repeated games with incomplete information or in splitting games, an important property is   that any martingale on a simplex  has bounded variations. This suggests the following.

\begin{defi}
A gambling house $\Gamma$ has vanishing $L_1$-variation if for every $\varepsilon>0$, there is $N$ such that for all $n\geq N$ and any sequence $(p_t)$ s.t. $p_{t+1} \in \widetilde{\Gamma}(p_t)$, one has $\frac{1}{n}\sum_{t=1}^n d_{KR}(p_{t+1},p_t) \leq \varepsilon $. 
\end{defi}

%The definition uses Cesaro-Mean, but one can extends the definition to Abel means or to any $\theta$ evaluations. 
 
The proof of the next result is inspired by Oliu-Barton \cite{Miquel} in the framework of splitting games. He shows that Player 1 (resp. Player 2) can uniformly guarantee any excessive-depressive function satisfying $MZ1$ (resp. $MZ2$). Our proof is much shorter because it uses the new $Q$-characterization.

\begin{pro} \label{thm5}
In a standard gambling game where $\Gamma$ has vanishing $L_1$-variation, if  $w$ in $\mathcal{B}(X\times Y)$    is excessive-depressive and satisfies $Q1$, then a strategy of Player 1 adapted to $w$   uniformly guarantees   $w$.
\end{pro}

\noindent \textbf{Proof:} Fix $\varepsilon>0$. 
Because Lipschitz functions are dense in the set of continuous functions,  there exists $K>0$ and a $K$-Lispchtiz function $u_{\varepsilon}$ that is uniformly $\varepsilon$-close to $u$. Consider  $w$ in $\mathcal{B}(X\times Y)$  be an excessive-depressive function satisfying $Q1$,   let $\sigma$ be a strategy of Player 1 adapted to $w$ and let  $\tau$ be any strategy of Player 2. We  fix the initial  states and write 
$\E=\E_{(x_1,y_1),\sigma,\tau}$.

 Then, the average payoff of the $n$-stage game is:

\begin{eqnarray*}
\frac{1}{n} \E \left(\sum_{t=1}^n  u(x_t,y_t)\right) &=& \frac{1}{n}  \E \left(\sum_{t=1}^n  u(x_t,y_t)-u(x_{t+1},y_t)\right)+\frac{1}{n}  \E \left(\sum_{t=1}^n  u(x_{t+1},y_t)\right)\\
&\geq &  \frac{1}{n}  \E \left(\sum_{t=1}^n  u_{\varepsilon}(x_t,y_t)-u_{\varepsilon}(x_{t+1},y_t)\right)-2\varepsilon+\frac{1}{n}  \E \left(\sum_{t=1}^n  u(x_{t+1},y_t)\right)\\
&\geq & - \frac{K}{n} \left(\sum_{t=1}^n  d_{KR}(p_{t+1},p_t)\right)-2\varepsilon+\frac{1}{n} \E \left(\sum_{t=1}^n  u(x_{t+1},y_t)\right)
\end{eqnarray*}

Because the gambling game is of vanishing variation, there is $N$ such that when $n\geq N$ one has $\frac{K}{n} \sum_{t=1}^n  d_{KR}(p_{t+1},p_t) \leq \varepsilon$. Because of $Q1$, $\E (u(x_{t+1},y_t)) \geq \E (w(x_{t+1},y_t)) = \E (w(x_{t},y_t))$.  Since  $w$ is depressive, $\E (w(x_{t},y_t) \geq \E (w(x_{t},y_{t-1}))$, so that $\E(w(x_{t+1},y_t)) \geq \E (w(x_{t},y_{t-1}))$, and this property holds for every $t$.  Consequently, $\E (u(x_{t+1},y_t)) \geq \E (w(x_{2},y_{1}))$ $= w(x_1,y_1)$. We obtain finally $\frac{1}{n} \E (\sum_{t=1}^n  u(x_t,y_t))\geq w(x_1,y_1) - 3 \varepsilon$,  ending the proof.

%Consequently,

\begin{cor} \label{cor17}
If $\Gamma$ and $\Lambda$ have vanishing $L_1$-variations then there is at most one excessive-depressive function  in $\mathcal{B}(X\times Y)$ satisfying $Q1$ and $Q2$. \end{cor}
\noindent \textbf{Proof:}  If such a function exists, it can be guaranteed by both players. So it must be the uniform value of the game, which is  unique whenever it exists. \\

%Thus, under vanishing  solution to the $Q$ characterization, 
Combining proposition \ref{pro2.5}, proposition \ref{thm5} and corollary \ref{cor17}, we obtain the existence of the uniform value in a class of gambling games:

\begin{thm} \label{corthm5}
In a standard and   idempotent gambling game where  $\Gamma$ and $\Lambda$ have vanishing $L_1$-variation,   the uniform value $v$ exists, and strategies adapted to $v$ are uniformly optimal. Moreover, $v$ is the unique excessive-depressive function   in $\mathcal{B}(X\times Y)$ satisfying $Q1$ and $Q2$. \end{thm}

\noindent \textbf{Proof:}  From proposition \ref{pro2.5}, any accumulation point of $v_{\lambda}$ satisfies $P1$ and $P2$, i.e.  $Q1$ and $Q2$,  so  is unique and is the uniform value, as shown above.\\
%Moreover, players, have stationary uniformly optimal strategies and $v_{\theta^n}$ uniformly converges to $v$ as $\theta^n_1 \rightarrow 0 $ for any sequence $\theta^n=(\theta^n_t)$ of decreasing evaluations (i.e. $\theta^n_{t+1} \leq \theta^n_t $ for every $t$). This holds true for any stochastic game because the uniform value exists, because for any decreasing evaluation $\theta$, $\sum_t \theta_t u(x_t,y_t)$ belongs to the convex envelop $\{\frac{1}{n} \sum_1^n u(x_t,y_t)\}_{n\in \mathbf{N}^*}$.

%Our proof is simpler in a sense because it is more easy to work with $P1$ than $MZ1$.

Observe that acyclicity is not assumed in theorem \ref{corthm5}. But vanishing $L_1$-variation is a form of acyclicity (it rules out for example non-constant periodic orbits). A more formal link between acyclicity and vanishing variation is given in  section \ref{sub47} of the Appendix.

\section{Open problems and future directions}

We introduce the class of gambling games. It is a sub-class of stochastic games which includes MDP problems, splitting games and product stochastic games. We define  a strong notion of acyclicity  under which we prove existence of the asymptotic value $v$ and we establish several characterizations of $v$ which are linked to the Mertens-Zamir system of functional equations (re-formulated in our more general set-up). We also prove that our condition is tight: a slight weakening of acyclicity implies non-existence of the asymptotic value. Our example is the first in the class of product stochastic games and is probably the simplest known counterexample of convergence for finite state spaces and compact action set (the first counterexample in this class was established by Vigeral \cite{Vigeral_2013}). Many questions merit to be investigated in a future research: 
\begin{itemize}
\item In standard gambling games, is it possible to characterize the asymptotic value in models where we know it exists (for example when $X$ and $Y$ are finite, transition function is polynomial and $\Gamma$ and $\Lambda$ are definable \cite{BGV_2015} in an o-minimal structure)? We know that the asymptotic value is excessive-depressive and satisfy P1 and P2, but this does not fully characterize it.   
\item Is there an asymptotic value if one house is strongly acyclic and the other not necessarily leavable? As seen, even when both houses are weakly acyclic, we may have divergence:  strong acyclicity of one of the two houses is necessary. 
\item Is there a strongly  acyclic gambling game where the uniform value fails to exist?
\item Recently splitting games have been extended to continuous time and linked to differential games with incomplete information \cite{GR_2017}.  How our model and results extend to continuous time? The one player game has been investigated \cite{BQR, LQR}.
  
%\item how to extend the model to the dependent case (as dependent splitting games in Oliu-Barton, or the dependent incomplete information game in Mertens and Zamir).
\end{itemize}

%As Bolte, Gaubert and Vigeral proved, when the state space is finite, definability of the Shapley operator guarantee existence of the asymptotic value.
 
%$\bullet$ Autre exemples ??: Casino  ou chacun veut être plus riche. 
%
%Courses (tennis, on peut faire tous les tournois pour gagner des points au risque de se blesser, ou seulement les gros tournois), 
%
% Entreprises qui veulent être première sur le marché, investir en augmentant les prix, ou gagner des parts de marché ?
% 
%$\bullet$  important: Valeur uniforme, stratégies optimales simples ($\varepsilon$-optimales quand $\lambda$ est petit). Un intérêt de la caractérisation est qu'elle met en évidence des bonnes stratégies (comme chez  Miquel)
%
%$\bullet$ voir si le théorème marche avec un des 2 joueurs strongly acyclique seulement. 
%
% 
% $\bullet$  problème ouvert:  dire à un moment que c'est beaucoup plus compliqué quand pas leavable (pour 1 joueur dès que non expansive, existence déduit de Renault 2011, caractérisation Renault Venel ), problème ouvert intéressant: acyclic gambling games non leavable, au sens ? : there exists $\varphi:X \to \R$ lower semi-continuous, and $\psi:Y \to \R$ upper semi-continuous   such that: $\forall x \in X, p\in \Gamma^{\infty}(x)$, if $p\neq \delta_x$ then $\varphi(p)<\varphi(x)$  ? ceci ne peut pas marcher je pense (on peut en discuter)
% 
% 
% 
%
%
%Citer: Janos, BZ, Guillaume V, Bolte etc.... Laraki, MZ
%
%
%
%\newpage

\section{Appendix}

 \subsection{Proof of proposition \ref{pro1}}
$u$ being uniformly continuous over the compact set $X\times Y$, we consider a concave modulus of continuity $\omega: \R_+ \longrightarrow \R_+$:
$$|u(x,y)-u(x',y')|\leq \omega (d(x,x')+d(y,y')),\; \forall x,   x' \in X,  \forall y,y' \in Y.$$
$\omega$  is non decreasing, concave and  $\lim_{0}\; \omega=0$.  Denote by ${\cal C}$ the set of functions $v$ in $\mathcal{C}(X\times Y)$ satisfying: $|v(x,y)-v(x',y')|\leq \omega (d(x,x')+d(y,y')),\; \forall x,   x' \in X,  \forall y,y' \in Y.$  We start with a lemma. 

 \begin{lem} \label{lem1}For $v$ in ${\cal C}$, $p$, $p'$ in $\Delta(X)$, $q$, $q'$ in $\Delta(Y)$, $$|\tilde{v}(p,q)-\tilde{v}(p',q')|\leq \omega (d_{KR}(p,p')+d_{KR}(q,q')).$$\end{lem}

\noindent{Proof of lemma \ref{lem1}: By the Kantorovich duality theorem, there exists $\mu$ in $\Delta(X\times X)$ with first marginal $p$ and second marginal $p'$ satisfying: $d_{KR}(p,p')=\int_{X\times X}d(x,x') d\mu(x,x').$ Similarly there exists $\nu$ in $\Delta(Y\times Y)$ with first marginal $q$ and second marginal $q'$ satisfying: $d_{KR}(q,q')=\int_{Y\times Y}d(y,y') d\nu(y,y').$ We have for all $x$, $x'$, $y$, $y'$:
$$v(x,y)\geq v(x',y')- \omega(d(x,x')+d(y,y')).$$
We  integrate the above inequality with respect to the  probability $\mu \otimes \nu$, and obtain using the concavity of $\omega$:
\begin{eqnarray*}
\tilde{v}(p,q) & \geq &\tilde{v}(p',q')- \int_{X^2 \times Y^2} \omega(d(x,x')+d(y,y'))\; d\mu(x,x')d\nu(y,y'),\\
 & \geq &\tilde{v}(p',q') - \omega\left( \int_{X^2 \times Y^2} d(x,x')+d(y,y') \; d\mu(x,x')d\nu(y,y')\right)\\
 &=&   \tilde{v}(p',q') - \omega\left(d_{KR}(p,p')+d_{KR}(q,q')\right).
 \end{eqnarray*}
 
 \vspace{0,5cm}
\noindent We now return to the proof of proposition \ref{pro1}. Fix $\lambda$ in $(0,1]$. Given $v$ in ${\cal C}$, define $\Phi(v):X \times Y \longrightarrow \R$ by:
$\Phi(v)(x,y)=\sup_{p \in \Gamma(x)}\inf_{q \in \Lambda(y)}  \lambda \; u(x,y)+ (1-\lambda) \; \tilde{v}(p,q).$

\noindent Consider the zero-sum game with strategy spaces $\Gamma(x)$ and $\Lambda(y)$ and payoff function  $(p,q)\mapsto  \lambda \; u(x,y)+ (1-\lambda) \; \tilde{v}(p,q)$. The strategy spaces  are convex compact and the payoff function is  is continuous and affine in each variable, hence by Sion's theorem  we have:
$$\Phi(v)(x,y)=\max_{p \in \Gamma(x)}\min_{q \in \Lambda(y)}  \lambda \; u(x,y)+ (1-\lambda) \; \tilde{v}(p,q)=\min_{q \in \Lambda(y)} \max_{p \in \Gamma(x)} \lambda \; u(x,y)+ (1-\lambda) \; \tilde{v}(p,q).$$

Consider $(x,y)$ and $(x',y')$ in $X\times Y$, and let $p$ in $\Gamma(x)$ be an optimal strategy of Player 1 in the zero-sum game corresponding to $(x,y)$. The gambling game has non expansive transitions, so there exists $p'\in \Gamma(x')$ such that $d_{KR}(p,p')\leq d(x,x')$. Consider any $q'$ in $\Lambda(y')$, there exists $q$ in $\Gamma(y)$ with  $d_{KR}(q,q')\leq d(y,y')$. Now, using lemma \ref{lem1} we write: 

$ \lambda u(x',y') +(1-\lambda) \tilde{v}(p',q')$

$\geq \lambda (u(x,y)-\omega(d(x,x')+d(y,y')) + (1-\lambda) ( \tilde{v}(p,q) -\omega({d_{KR}(p,p')+d_{KR}(q,q')})),$

$ \geq   \Phi(v)(x,y)-\omega({d(x,x')+d(y,y')}).$
 
 We obtain $\Phi(v)(x',y')\geq  \Phi(v)(x,y)-\omega({d(x,x')+d(y,y'))}$, and   $\Phi(v)$ belongs to ${\cal C}$.\\
 
 The rest of the proof is very standard. ${\cal C}$ is a complete metric space for $\|v-w\|=\sup_{(x,y)\in X\times Y} |v(x,y)-w(x,y)|$, and we have  $\|\Phi(v)-\Phi(w)\|\leq (1-\lambda) \|v-w\|$, so $\Phi$ is $(1-\lambda)$-contracting. Hence $\Phi$ has a unique fixed point which is $v_\lambda$. Each $v_\lambda$ is in ${\cal C}$, and we obtain that the family $(v_\lambda)_{\lambda\in (0,1]}$ is equicontinuous, ending the proof of proposition \ref{pro1}.

 \subsection{Proof of proposition \ref{pro11}} We proceed in 4 steps.\\
  
\noindent{\bf 1. The game at $(b,a')$:} It is intuitively clear that $y_\lambda \geq x_\lambda$ since Player 1 is better off when the players have different locations. We now formalize this idea. Consider the game at $(b,a')$. The current payoff is 1, and Player 1 has the option not to move, so we obtain by definition \ref{def1}:
 $$y_\lambda\geq \lambda+(1-\lambda) \min_{\beta \in J}(\beta x_\lambda +(1-\beta-\beta^2)y_\lambda +\beta^2 1).$$
 \noindent This implies $$\lambda y_\lambda \geq \lambda +(1-\lambda)  \min_{\beta \in J}(\beta (x_\lambda -y_\lambda) +\beta^2 (1-y_\lambda)),$$ and since $y_\lambda\leq 1$, we obtain: $x_\lambda\leq y_\lambda$. Now,  $\min_{\beta \in J}(\beta (x_\lambda -y_\lambda) +\beta^2 (1-y_\lambda))\geq \min_{\beta \in J}\beta (x_\lambda -y_\lambda) +\min_{\beta \in J}\beta^2 (1-y_\lambda)=1/4(x_\lambda-y_\lambda)$, hence:
    \begin{equation} \label{eq11}  (1-\lambda)(y_\lambda-x_\lambda)  \geq 4 \lambda (1-y_\lambda).\end{equation}
   In the same spirit, in the game at $(a,a')$, Player 2  has the option not to move, so we have:$$x_\lambda\leq (1-\lambda) \max_{\alpha \in I} (\alpha y_\lambda +(1-\alpha-\alpha^2)x_\lambda + \alpha^2 z_\lambda).$$ Hence, 
\begin{equation} \label{eq22}\lambda x_\lambda \leq (1-\lambda) \max_{ \alpha \in I} (\alpha(y_\lambda-x_\lambda) + \alpha^2(z_\lambda-x_\lambda)).\end{equation}

 \noindent{\bf 2. The game at $(a,a')$:} Consider now the game at $(a,a')$.  By definition \ref{def1}, $x_\lambda$ is the value of the game (possibly played with mixed strategies), where Player 1 chooses $\alpha$ in $I$, Player 2 chooses $\beta$ in $J$ and the payoff to Player 1 is: $(1-\lambda) g_\lambda(\alpha,\beta)$, where $g_\lambda(\alpha,\beta)=$
 $$x_\lambda((1-\alpha-\alpha^2)(1-\beta-\beta^2)+\alpha \beta) + y_\lambda (\beta(1-\alpha-\alpha^2)+\alpha(1-\beta-\beta^2)) + \beta^2 (1-\alpha^2) 1 + \alpha^2(1-\beta^2) z_\lambda.$$
\noindent We want to prove that in this game, it is a dominant strategy  for Player 2 not to move, that is to choose $\beta=0$. We need to show that for all $\alpha$ and $\beta$, $g_\lambda(\alpha, 0)\leq g_\lambda(\alpha, \beta)$. As a function of $\beta$, $g_\lambda(\alpha,\beta)$ can be written as a constant plus:
$$\beta(1-2\alpha-\alpha^2)(y_\lambda-x_\lambda) + \beta^2 (-(1-\alpha-\alpha^2)x_\lambda- \alpha y_\lambda+1-\alpha^2-\alpha^2 z_\lambda).$$
So we want to show that for all $\alpha$ in $I$, $ \beta$ in $J$:
$$ (1-2\alpha-\alpha^2-\alpha \beta)(y_\lambda-x_\lambda) + \beta  ((1- \alpha^2)(1-x_\lambda)- \alpha^2 z_\lambda)\geq 0$$
Since the expression is decreasing in $\alpha$, it is enough to prove it with $\alpha=1/4$:
$$ (7-4 \beta)(y_\lambda-x_\lambda) + \beta  (15(1-x_\lambda)- \ z_\lambda)\geq 0.$$
This is true for $\beta=0$, and will be true for all $\beta$ in $J$ if and only if it is true for $\beta=1/4$, so we are left with proving:
\begin{equation}\label{eq21}  24(y_\lambda-x_\lambda) +   15(1-x_\lambda)-  \ z_\lambda\geq 0.\end{equation}

Consider $\lambda\leq 1/32$, and recall that $z_\lambda\leq 16 \lambda$.  If $x_\lambda\leq 1/2$, then clearly (\ref{eq21}) holds. Assume on the contrary that $x_\lambda\geq 1/2$, then   $z_\lambda \leq x_\lambda$, and (\ref{eq22}) gives: $\lambda x_\lambda \leq (1-\lambda)\frac{1}{4} (y_\lambda-x_\lambda) + 0$, so $(y_\lambda-x_\lambda)\geq 2 \lambda$ and (\ref{eq21}) holds as well. 

We have shown that in the $\lambda$-discounted game at $(a,a')$ with $\lambda\leq 1/32$,  Player 2 has a  pure dominant strategy which is $\beta=0$. Considering a pure best reply of Player 1 against this strategy  implies that the game at $(a,a')$ has a value in pure strategies satisfying $x_\lambda=  (1-\lambda) \max_{\alpha \in I} (\alpha y_\lambda +(1-\alpha-\alpha^2)x_\lambda + \alpha^2 z_\lambda)$, i.e. equation (\ref{eq32}) is proved. \\

   \noindent{\bf 3. Small discount factors:}
   For the sake of contradiction, assume  that $z_{\lambda_n}\geq x_{\lambda_n}$ for a vanishing sequence $\lambda_n$ of discount factors. Then by equation (\ref{eq32}), we have for each $n$: 
   $$ \lambda_n  x_{\lambda_n} = (1-\lambda_n)  (\frac{1}{4}(y_{\lambda_n}-x_{\lambda_n}) + \frac{1}{16}(z_{\lambda_n}-x_{\lambda_n})).$$
  Since $z_{\lambda_n}$ converges to 0, so does   $x_{\lambda_n}$ and $y_{\lambda_n}$, and moreover $\frac{y_{\lambda_n}-x_{\lambda_n}}{\lambda_n}$ converges to 0. This is in contradiction with equation (\ref{eq11}). We have shown (\ref{eq31}).\\

    \noindent{\bf 4. The game at $(b,a')$ again:} We  proceed as for the game at $(a,a')$ and will show that in the game at $(b,a')$, it is a dominant strategy for Player 1 not to move. By definition, $y_\lambda$ is the value of the game where Player 1 chooses $\alpha$ in $I$, Player 2 chooses $\beta$ in $J$ and the payoff is $\lambda +(1-\lambda)h_\lambda(\alpha, \beta)$, with $h_\lambda(\alpha, \beta)=$
     $$y_\lambda((1-\alpha-\alpha^2)(1-\beta-\beta^2)+\alpha \beta) + x_\lambda (\beta(1-\alpha-\alpha^2)+\alpha(1-\beta-\beta^2)) + \beta^2 (1-\alpha^2) 1 + \alpha^2(1-\beta^2) z_\lambda.$$
We want to show that     $h_\lambda(0,\beta)\geq h_\lambda(\alpha, \beta)$ for all $\alpha$ and $\beta$. That is, for all $\alpha$ and $\beta$,
$$(x_\lambda-y_\lambda)(1-2\beta-\alpha \beta-\beta^2) + \alpha((1-\beta^2)(z_\lambda-y_\lambda)-\beta^2) \leq 0.$$
    For $\lambda$ small enough, we have  $z_\lambda \leq x_\lambda\leq y_\lambda$, and the above    property is satisfied. Hence in the game at $(b,a')$ it is dominant  for Player 1 to choose $\alpha=0$. Consequently, Player 2 has a pure optimal strategy and we can write:
 $$y_\lambda= \lambda+(1-\lambda) \min_{\beta \in J}(\beta x_\lambda +(1-\beta-\beta^2)y_\lambda +\beta^2),$$
proving  equation (\ref{eq33}). And  the proof of proposition \ref{pro11} is complete.

\subsection{$L_2$-variation, $L_1$-variation and acyclicity} \label{sub47}

\begin{defi}
A gambling house $\Gamma$ is of bounded $L_2$-variation if there is $C>0$ such that for every sequences $\{p_t\}$ satisfying $p_{t+1} \in \widetilde{\Gamma}(p_t)$, one has $$\sum_{t=1}^{\infty} d_{KL}(p_{t+1},p_t)^2 \leq C <+\infty $$
\end{defi}

For example splitting games   have bounded $L_2$-variation.

\begin{pro}
If $\Gamma$ is idempotent, non-expansive, leavable and of bounded $L_2$-variation, then it is acyclic and has vanishing $L_1$-variation.
\end{pro}

\noindent \textbf{Proof:} Bounded $L_2$-variation $\implies$ weak acyclicity because the real valued function: $$\varphi(x)=\sup_{\{(p_t) \ \textit{s.t.} \ p_0=x \ \textit{and} \ p_{t+1} \in \widetilde{\Gamma} (p_t)\}}  \sum_{t=1}^{\infty} d_{KL}(p_{t+1},p_t)^2,$$
is strictly decreasing along non-constant orbits of $\widetilde{\Gamma}$ (i.e. $\arg \max_{p\in \Gamma(x)} \varphi (p)= \delta_x$).  But for idempotent $\Gamma$, weak acyclicity and acyclicity coincide because $\Gamma^{\infty}=\Gamma$. Continuity of $\varphi$ is a consequence of non-expansivity of $\Gamma$ and the fact the   bound $C$ is uniform over the sequences $\{p_t\}$ satisfying $p_{t+1} \in \widetilde{\Gamma}(p_t)$.   Finally, that bounded $L_2$-variation implies vanishing $L_1$-variation is a consequence of Cauchy-Schwartz inequality ($\frac{1}{n}\sum_{t=1}^n d_{KL}(p_{t+1},p_t) \leq \frac{1}{\sqrt{n}} \sqrt{\sum_{t=1}^{n} d_{KL}(p_{t+1},p_t)^2}$).

%Also, it is well known that when the uniform value $v$ exists, it is the asymptotic value for all decreasing evaluations (see for instance Sorin 2002).


\begin{thebibliography}{99}

\bibitem{AB_2006} Aliprantis C. D.  and K.C. Border. \textit{Infinite Dimensional Analysis: a Hitchhiker's Guide}.\newblock Springer (2006)


\bibitem{AM_1995}
Aumann R., M. Maschler and R. Stearns. \textit{Repeated Games with Incomplete
  Information}.\newblock The MIT Press (1995).
  
  \bibitem{BK} Bewley T. and E. Kohlberg. The Asymptotic Theory of
Stochastic Games. \textit{Mathematics of Operations Research}, \textbf{1},
197-208 (1976).


   \bibitem{BGV_2015} Bolte J., S. Gaubert and G. Vigeral. \newblock Definable Zero-Sum Stochastic Games. \textit{Mathematics of Operations Research}, \textbf{40}, 171-191 (2015).

   \bibitem{BQR} Buckdahn R., M. Quincampoix and J. Renault. On Representation Formulas for Long Run Averaging Optimal Control Problem. \textit{Journal of Differential Equations} \textbf{259}, 5554-5581 (2015).

  
   \bibitem{DS_1965} Dubins L. E. and L. J. Savage. \textit{Inequalities for Stochastic
Processes}, McGraw-Hill. 2nd edition 1976, Dover  (1965).

   \bibitem{DMS} Dubins L. E., A. P. Maitra and W. D. Sudderth. Invariant Gambling Problems and Markov Decision Processes.
E. A. Feinberg, A. Shwartz, eds. \textit{Handbook of Markov Decision Processes}. Kluwer, 409-428 (2002).

\bibitem{Everett} Everett H. Recursive Games. \textit{Contributions to the
Theory of Games, III}, Dresher M., A.W. Tucker and P. Wolfe (eds.), Annals
of Mathematical Studies, \textbf{39}, Princeton University Press, 47-78 (1957).


\bibitem{FSV_2008}Flesch J., G. Schoenmakers and K. Vrieze.  Stochastic Games on a Product State Space. \textit{Mathematics of Operations Research}
\textbf{33}(2), 403-420 (2008).

\bibitem{GR_2017} Gensbittel F. and C.Rainer. A Probabilistic Representation for the Value of Zero-Sum Differential Games with Incomplete Information on Both Sides.
to appear in \textit{SIAM J. Control Optim} (2017).

\bibitem{GR_2015} Gensbittel F. and J. Renault. The value of Markov chain games with lack of information on both sides. \textit{Mathematics of Operations Research} \textbf{40}(4), 820--841(2015). 


\bibitem{Kohlberg_1974} Kohlberg E. Repeated Games with Absorbing States. \textit{%
Annals of Statistics}, \textbf{2}, 724-738 (1974).

\bibitem{KG}  Kamenica E.  and M. Gentzkow. Bayesian Persuasion. \textit{The American Economic Review}, \textbf{101(6)}, 2590-2615 (2011).


\bibitem{Laraki_2001a} Laraki R. Variational Inequalities, Systems of
Functional Equations and Incomplete Information Repeated Games. \textit{SIAM
J. Control and Optimization}, \textbf{40}, 516-524 (2001).


\bibitem{Laraki_2001b} Laraki R. The Splitting Game and Applications.
\textit{International Journal of Game Theory}, \textbf{30}, 359-376.



\bibitem{LS_2014} Laraki R. and S. Sorin. Advances in Zero-Sum Dynamic Games. Chapter 2 in \textit{Handbook of Game Theory IV}, edited by H. P. Young and S. Zamir, 27-93 (2014).


\bibitem{LS_2004} Laraki R. and W. D. Sudderth. The Preservation of Continuity and Lipschitz Continuity by Optimal Rewards Operators. \textit{Mathematics of Operations Research}, \textbf{29}, 672-685 (2004).

\bibitem{LQR}  LI X., M. Quincampoix and J. Renault. Limit Value for Optimal Control with General Means.  \textit{Discrete and Continuous Dynamical System A}, \textbf{36}(4), 2113-2132 (2016).

\bibitem{MS_1996}
Maitra A. and W. D. Sudderth. \textit{Discrete Gambling and Stochastic Games}, \textbf{32}. \newblock Springer Verlag (1996).

\bibitem{Mertens_1987}
Mertens J.-F.. Repeated Games. In \textit{Proceedings of the International
  Congress of Mathematicians}, {V}ol. 1, 2 ({B}erkeley, {C}alif., 1986),
  (Providence, RI), pp.~1528--1577, Amer. Math. Soc. (1987).

\bibitem{MNR_2009} Mertens J.-F., A. Neyman and D. Rosenberg. Absorbing
Games with Compact Action Spaces, \textit{Mathematics of Operations Research}%
, \textbf{34}, 257-262 (2009).

\bibitem{MertensSorinZamir_2015}
Mertens J.-F., S. Sorin  and S. Zamir. \textit{Repeated Games}. Cambridge University Press (2015).
  
\bibitem{MertensZamir_1971} Mertens J.-F. and S. Zamir. The Value of Two-Person
Zero-Sum Repeated Games with Lack of Information on Both Sides. \textit{%
International Journal of Game Theory}, \textbf{1}, 39-64 (1971).

\bibitem{NeymanSorin_2003} Neyman A. and S. Sorin (eds.)  \textit{Stochastic
Games and Applications}, NATO Science Series C 570, Kluwer Academic
Publishers (2003).

\bibitem{Miquel}Oliu Barton M. The Splitting Game: Uniform Value and Optimal Strategies. To appear in \textit{Dynamic Games and Applications} (2017)

\bibitem{Renault_2011} Renault J. Uniform Value in Dynamic Programming. \textit{Journal of the European Mathematical Society} \textbf{13}(2), 309-330 (2011)

\bibitem{RV_2013}Renault J. and X. Venel. A Distance for Probability Spaces, and Long-Term Values in Markov Decision Processes and Repeated Games. arXiv:1202.6259. To appear in M\textit{athematics of Operations Research} (published online November 28, 2016). 
  
\bibitem{RS_2001} Rosenberg D. and S. Sorin. An Operator Approach to Zero-Sum Repeated Games. \textit{Israel Journal of Mathematics}, \textbf{%
121}, 221- 246 (2001) .

\bibitem{Sh} Shapley L. S. Stochastic Dames. \textit{Proceedings of
the National Academy of Sciences of the U.S.A}, \textbf{39}, 1095-1100 (1953).


\bibitem{Schal} Schal M. On Stochastic Dynamic Programming: A Bridge Between Markov Decision Processes and Gambling.
Markov Processes and Control Theory, \textit{Math. Res.} \textbf{54}, Akademie-Verlag, Berlin, Germany, 178-216 (1989).

\bibitem{Sorin_2002} Sorin S. \textit{A First Course on Zero-Sum Repeated Games}. \newblock Springer (2002).

 \bibitem{SorinVigeral_2013} Sorin S. and G. Vigeral. Existence of the Limit Value of Two Person Zeo-Sum Discounted Repeated Games via Comparison Theorems. \textit{JOTA} \textbf{157} (2), 564-576 (2013).

%\bibitem{SorinVigeral_2015} Sorin S. and G. Vigeral. Reversibility and Oscillations in Zero-Zum Discounted Stochastic Games. \textit{Journal of Dynamics and Games}, \textbf{2}(1), 103-115 (2015).
 
\bibitem{Vigeral_2013}Vigeral G. A Zero-Sum Stochastic Game with Compact Action Sets and no Asymptotic Value. \textit{Dynamic Games and Applications}, \textbf{3} (2), 172-186 (2013).

\bibitem{Ziliotto_2016} Ziliotto B. Zero-Sum Repeated Games: Counterexamples to the Existence of the Asymptotic Value and the Conjecture $maxmin=\lim v(n)$.
\textit{Annals of Probability}, \textbf{44}(2), 1107-1133 (2016).

\bibitem{Ziliotto_2015} Ziliotto B. A Tauberian Theorem for Nonexpansive Operators and Applications to Zero-Sum Stochastic Games. arXiv:1501.06525v2  2015. To appear in \textit{Mathematics of Operations Research}.


\end{thebibliography}
 \end{document}